\documentclass[10pt,a4paper,reqno]{amsart}
\usepackage{amsmath}
\usepackage{amsfonts}
\usepackage{amsthm}
\usepackage{mathtools}
\usepackage{amscd}
\usepackage[utf8]{inputenc}
\usepackage{t1enc}
\usepackage[mathscr]{eucal}
\usepackage{indentfirst}
\usepackage{graphicx}
\usepackage{graphics}
\usepackage{pict2e}
\usepackage{tikz-cd}
\usepackage{enumitem}
\usepackage{epic}
\usepackage{float}
\usepackage{MnSymbol}
\usepackage{multirow}
\usepackage{shuffle} 
\usepackage[dvipsnames]{xcolor}
\usepackage{hyperref}
\usepackage{comment}
\usepackage{tcolorbox}
\numberwithin{equation}{section}
\usepackage[margin=2.9cm]{geometry}
\usepackage{epstopdf} 
\usepackage[nocompress]{cite}

\allowdisplaybreaks

\newmuskip\pFqmuskip
\newcommand*\pFq[6][8]{%
	\begingroup 
	\pFqmuskip=#1mu\relax
	\mathchardef\normalcomma=\mathcode`,
	\mathcode`\,=\string"8000
	\begingroup\lccode`\~=`\,
	\lowercase{\endgroup\let~}\pFqcomma
	{}_{#2}F_{#3}{\left[\genfrac..{0pt}{}{#4}{#5};#6\right]}%
	\endgroup
}
\newcommand{\pFqcomma}{{\normalcomma}\mskip\pFqmuskip}

\theoremstyle{plain}
\newtheorem{theorem}{Theorem}[section]

\newtheorem{proposition}[theorem]{Proposition}

\theoremstyle{definition}
\newtheorem{definition}[theorem]{Definition}
\newtheorem{conjecture}[theorem]{Conjecture}
\newtheorem{remark}[theorem]{Remark}
\newtheorem{problem}[theorem]{Problem}

\newcommand{\MZV}[2]{\textsf{MZV}^{#1}_{#2}}
\DeclareMathOperator{\Li}{Li}

\renewcommand{\Im}{\operatorname{Im}}

\newtheorem{romexample}{Example}

\begin{document}
	
	\title{Discovering hypergeometric series with harmonic numbers via Wilf-Zeilberger seeds}
	
	\author[Kam Cheong Au]{Kam Cheong Au}
	
	\address{University of Cologne \\ Department of Mathematics and Computer Science \\ Weyertal 86-90, 50931 Cologne, Germany} 

\email{kau@uni-koeln.de}
	

	\begin{abstract}By extracting coefficients from Wilf-Zeilberger pairs with respect to auxiliary parameters, we discover many nontrivial hypergeometric series involving harmonic numbers. In particular, we obtain a rapidly convergent series for the depth-two multiple zeta value $\zeta(5,3)$, which appears to be the first result of its kind in the literature. We also investigate the Hilbert series attached to a WZ-seed and conjecture that it admits a remarkably simple form, suggesting the presence of an underlying graded-algebra structure associated with WZ-seeds.
	\end{abstract}
	
	\maketitle
	\section{Introduction}
	A \textit{Wilf-Zeilberger (WZ) pair} \cite{AequalsB, wilf1992algorithmic, wilf1990rational} consists of two functions $F$ and $G$ of two variables satisfying
\begin{equation}\label{WZ}
	F(n+1,k) - F(n,k) = G(n,k+1) - G(n,k).
\end{equation}
Originally introduced in the context of combinatorial summation identities, the WZ method has since found broad applications in the discovery and proof of hypergeometric series that are not readily accessible by classical techniques \cite{guillera2008hypergeometric, guillera2018dougall, guillera2003new, sun2024new}. A prime example is the rapidly convergent series for Apéry’s constant \cite{amdeberhan1998hypergeometric}:
\begin{equation}\label{aux_0}
	\sum_{n\geq 1} \left(\frac{-1}{2^{10}}\right)^n \frac{(1)_n^5}{(\frac12)_n^5} \frac{32-160 n+205 n^2}{n^5} = -2\zeta(3),
\end{equation}
where $(a)_n := a(a+1)\cdots (a+n-1)$ denotes the rising factorial.\par

The theory of WZ-pairs has recently advanced substantially with the author's introduction of the notion of a \textit{WZ-seed} \cite{au2025wilf, au2024wilfQ, au2025multiple}. This concept, which had already appeared implicitly (though unnamed) in earlier work of Gessel \cite{gessel1995finding}, allows for the systematic construction of WZ-pairs. This framework also unifies almost all WZ-style proofs of hypergeometric identities in the literature; see \cite{chu2014accelerating, chu2023pi, guillera2018dougall, hou2023taylor, wei2023some, mohammed2005infinite, pilehrood2011bivariate, pilehrood2010series, zhang2015common, chu2011dougall, pilehrood2008simultaneous}. One highlight of the application of WZ-seeds is their use in proving previously conjectural evaluations such as
\begin{align}\label{aux_7}\sum_{n\geq 1} \left(-\frac{2^{10}}{5^5}\right)^n \frac{(1)_n^9}{(\frac12)_n^5 (\frac15)_n(\frac25)_n(\frac35)_n(\frac45)_n} \frac{30-425 n+2275 n^2-5600 n^3+5532 n^4 }{n^9} &= -380928 \zeta(5), \\
	\sum_{n\geq 0} \left(\frac{1}{2^6}\right)^n \frac{(\frac12)_n^7}{(1)_n^7} (168 n^3+76 n^2+14 n+1) &= \frac{32}{\pi^3},
\end{align}
which were conjectured by Zhao\footnote{\url{https://mathoverflow.net/questions/281009}} and Gourevich \cite{cohen2021rational}, respectively. \par

WZ-seeds come with auxiliary parameters, denoted throughout this article by $a,b,c,d,\ldots$. These parameters enable us to apply coefficient extraction and thereby produce a rich collection of evaluations involving harmonic numbers. This is a central theme of the article.\par

The main text of the article consists of ten examples, each concerning a family of hypergeometric series with a fixed hypergeometric part. Before introducing our methodology and the computational details, we summarize some notable features revealed by our computations.

\subsection{Rapidly convergent series and multiple zeta values}
One immediate application, although of less theoretical interest, is the derivation of rapidly convergent series for zeta values such as $\zeta(5)$ and the multiple zeta value (MZV) $\zeta(5,3)$. Throughout, we use the convention
$$\zeta(s_1,\ldots,s_k) := \sum_{n_1 > \cdots > n_k \geq 1} \frac{1}{n_1^{s_1}\cdots n_k^{s_k}}.$$
We denote by $\MZV{}{n}$ the $\mathbb{Q}$-span of all MZVs of weight $n$. Among the many identities proved in this article, the following are particularly notable (Examples \ref{Ex_3} and \ref{Ex_4}):
\begin{align}\label{aux_00}&\sum_{n\geq 1}  \left(\frac{-1}{27}\right)^n \frac{(1)_n^2}{(\frac13)_n (\frac23)_n} \frac{17-136 n+408 n^2-640 n^3+560 n^4}{n^5 (2n-1)^4} = 180 \zeta (5)-\frac{56 \pi ^2 \zeta (3)}{3}, \\ &\sum_{n\geq 1} \left(\frac{-1}{2^{10}}\right)^n \frac{(1)_n^5}{(\frac12)_n^5} \left(\frac{351-1755 n+2240 n^2}{n^{10}}-\frac{\left(11 H_n^{(5)}-H_{2 n}^{(5)}\right) (205 n^2-160 n+32)}{n^5}\right) = -4 \zeta(5,3)-\frac{\pi ^8}{4725}. \nonumber\end{align}
They enable high-precision evaluations of $\zeta(5)$ and $\zeta(5,3)$, respectively. The first identity was conjectured by Sun.\footnote{\url{https://mathoverflow.net/questions/486353}} The rapidly convergent series for $\zeta(5,3)$, a genuine depth-two MZV, appears to be the first of its kind in the literature. \par
WZ-seeds can systematically generate many more complicated rapidly convergent series for certain mathematical constants; see \cite{zuniga2026fast, zuniga2025fast1, zuniga2025fast2}. 

It is natural to ask why higher-depth MZVs appear in such WZ-style identities. Their appearance here can be traced to coefficient extraction from a very-well-poised $_7F_6$ hypergeometric series. Specifically, consider the hypergeometric series from Examples \ref{Ex_1} and \ref{Ex_3}:
$$f(a,b,c,d,e) := \sum_{k\geq 0} \frac{(a+2 e+2 k+2) (a+e+1)_k (b+e+1)_k (c+e+1)_k (d+e+1)_k}{(e+1)_{k+1}(a-b+e+1)_{k+1} (a-c+e+1)_{k+1} (a-d+e+1)_{k+1}},$$
which defines an analytic function in a neighborhood of $(a,b,c,d,e) = (0,0,0,0,0)$. If $c(i_1,\ldots,i_5)$ denotes the corresponding power-series coefficient, so that
$$f(a,b,c,d,e) = \sum_{i_1,\ldots,i_5\geq 0} c(i_1,\ldots,i_5) a^{i_1} b^{i_2} \cdots e^{i_5},$$
then it is readily seen that $c(i_1,\ldots,i_5) \in \MZV{}{2+i_1+\cdots+i_5}$. This leads us to the following conjecture.
\begin{conjecture}
		Let $$V_n := \operatorname{Span}_\mathbb{Q}\{c(i_1,\ldots,i_5) \mid i_1+\cdots+i_5 = n\}.$$
		Then its Hilbert--Poincaré series is a rational function; more precisely,
	$$\sum_{n=0}^\infty \dim_\mathbb{Q} V_n t^n \stackrel{?}{\in} \mathbb{Q}(t).$$
\end{conjecture}
As is usual for conjectures of this type, it may be more natural to investigate a motivic counterpart, which would avoid intractable questions of linear independence \cite{gil2017multiple}. 
\subsection{Collision phenomenon from different WZ-pairs}
The fact that multiple WZ-pairs can yield the same hypergeometric part is an interesting phenomenon: together, they produce additional linear relations and can prove deeper hypergeometric identities than any single WZ-pair can. For example, in Example \ref{Ex_4}, we will see that all three WZ-pairs\footnote{Here we use explicit WZ-seeds from Section~2.}
$$\begin{aligned}
	F_1(n,k) &= \textsf{Dougall5F4}\left(a+n+\frac{1}{2},b+\frac{1}{2},c,d-n,k+d+n+\frac{1}{2}\right),\\
	F_2(n,k) &= \textsf{Dougall5F4}\left(a,b,c,d-n,e+k+2n+1\right),\\
	F_3(n,k) &= \textsf{Seed3}\left(a,b-n,c,k+d+n+1\right),
\end{aligned}$$
give hypergeometric identities with the configuration $a_n := \left(\frac{-1}{27}\right)^n \frac{(1)_n^2}{(\frac13)_n (\frac23)_n}$. In particular, for monomials of degree $4$, we obtain various identities of the form
$$\sum_{n\geq 1} a_n \left(\text{some rational terms involving harmonic numbers}\right) \in \MZV{}{7}.$$
The first WZ-pair, $F_1$, gives $46$ linearly independent identities of the above type; combining $F_1$ and $F_2$ gives $54$ such identities; and using all three WZ-pairs increases this number to $59$; see Table \ref{Table_Ex_4}. For instance, one such identity is
$$\begin{aligned}
	&\sum_{n\geq 1} a_n\left( \frac{9 (-1+3 n)}{n^6 (-1+2 n)^2}-\frac{(5-32 n+56 n^2) (H_n^{(4)}-4 H_{2 n}^{(4)})}{n^3 (-1+2 n)^2}\right) = -\frac{\pi ^4 \zeta (3)}{5}+\frac{62 \pi ^2 \zeta (5)}{3}-205 \zeta (7).
\end{aligned}$$
The identity for $\zeta(5)$ in \eqref{aux_00} cannot be derived by considering $F_1$ alone. \par

For most hypergeometric configurations, such collisions do not occur. For example, $\left(\frac{1}{2^6}\right)^n \frac{(\frac12)_n^7}{(1)_n^7}$ in \eqref{aux_7} can be produced by only a single WZ-pair. The configurations \footnote{or its up-side-down version, see \cite{guillera2014ramanujan}} for which collisions do occur, such as $\left(\frac{-1}{27}\right)^n \frac{(1)_n^2}{(\frac13)_n (\frac23)_n}$ above and $\left(\frac{1}{2^6}\right)^n \frac{(\frac12)_n^3}{(1)_n^3}$ in Example \ref{Ex_7}, seem invariably to be connected to modular forms. This suggests that the corresponding identities involving harmonic numbers have a much richer structure.

\subsection{Relation to hypergeometric modular parameterization}
The hypergeometric series
$$\sum_{n\geq 0} \frac{(\frac12)_n^3}{(1)_n^3} z^n$$
is one of the few known hypergeometric series that admit a modular parameterization \cite{borwein1987pi, allen2025explicit}. Specialization at a CM point gives $z=1/2^6$. This particular value, $z=1/64$, is noteworthy because it lies at the intersection of the WZ approach used in this article and the modular approach.
The base case is the well-known Ramanujan-type $1/\pi$ formula
$$\sum_{n\geq 0} \left(\frac{1}{2^6}\right)^n \frac{(\frac12)_n^3}{(1)_n^3} (42n+5) = \frac{16}{\pi},$$
for which both the WZ proof \cite{guillera2008hypergeometric} and the modular proof are standard. \par
Here the collision phenomenon occurs again: three different WZ-seeds generate three WZ-pairs that provide linear relations for series with the hypergeometric part $a_n := \left(\frac{1}{2^6}\right)^n \frac{(\frac12)_n^3}{(1)_n^3}$; see Example \ref{Ex_7} below. Some simpler instances involving harmonic numbers, such as
$$\begin{aligned}
	&\sum_{n\geq 0} a_n\left(\frac{31+162 n}{3+18 n}+\left(-\frac{H_n}{2}-H_{3 n}+2 H_{6 n}\right) (5+42 n)\right) \stackrel{?}{=} \frac{48 \log (2)}{\pi }
\end{aligned}$$
can be proved by both the WZ and modular approaches \cite{zhou2025notes}. \par
Other identities, such as
$$\begin{aligned}
	&\sum_{n\geq 0} a_n \left(-\frac{9}{(1+2 n)^2}+\frac{1}{3} \left(-2 H_n^{(3)}+17 H_{2 n}^{(3)}\right) (42n+5)\right) = \frac{80 \zeta (3)}{\pi }-\frac{128 G}{3}\\
	&\sum_{n\geq 0} a_n \left( \frac{5 \left(-673+270 (1+2 n)^2 H_{2 n}^{(2)}\right)}{27 (1+2 n)^3}+\left(\frac{H_n^{(4)}}{16}+\left(H_{2 n}^{(2)}\right){}^2\right) (5+42 n) \right) = -\frac{3251 \pi ^3}{810}
\end{aligned}$$
currently have no accessible proof using modular forms and were discovered only through the WZ strategy; here $G$ denotes Catalan's constant. Still other identities exhibit the converse phenomenon, such as the following identity conjectured by Sun \cite{sun2022conjectures}:
$$\sum_{n\geq 0}  a_n (42n+5) \left(H_{2n}^{(3)} - \frac{43}{352} H_n^{(3)}\right) =\frac{555\zeta(3)}{77\pi} - \frac{32G}{11}.$$
It remains inaccessible to our WZ approach despite the availability of three WZ-pairs that yield linear relations. Recently, Zhou obtained a modular proof of this identity in \cite{sun2026series}. \par
Based on these examples, we raise the following question as a potentially interesting direction for further investigation.

\begin{problem}
What is the relation between the WZ and modular approaches to proving hypergeometric identities? How can the two approaches complement each other in proving identities that are inaccessible to either approach alone?
\end{problem}

\subsection{Hilbert series and the ``symmetry group'' of a WZ-seed}
WZ-seeds also seem to possess certain symmetries that manifest themselves in the corresponding infinite series. Let us revisit the identity \eqref{aux_0} and set $a_n := \left(\frac{-1}{2^{10}}\right)^n \frac{(1)_n^5}{(\frac12)_n^5}$. This identity comes from a WZ-seed with four auxiliary parameters, so we have five parameters with respect to which coefficients can be extracted (Example \ref{Ex_3}). For a monomial of degree $d$, one obtains identities of the form
\begin{equation}\label{aux_15}\sum_{n\geq 1} a_n \left(\text{some rational terms involving harmonic numbers}\right) \in \MZV{}{3+d}.\end{equation}
Our computations suggest the following statement (Conjecture \ref{gen_conj_Ex3}).
\begin{conjecture}Let $\dim V_d(S)$ denote the number of linearly independent identities\footnote{See the statement of the conjecture in the main text for a rigorous formulation.} obtained by comparing monomials of degree $d$ as in \eqref{aux_15}. Then
$$\sum_{d\geq 0} \dim V_d(S) t^d \stackrel{?}{=} \frac{1}{(1-t)(1-t^2)(1-t^3)(1-t^4)(1-t^5)} = 1+t+2 t^2+3 t^3 +\cdots.$$
\end{conjecture}
This suggests that the space of possible identities obtainable from this five-parameter WZ-pair is, in some sense, isomorphic to $\mathbb{Q}[e_1,\ldots,e_5]$, with $e_i$ assigned degree $i$. Indeed, our computations also reveal analogues of these putative ``generators'' among the infinite-series identities: identity~\eqref{aux_0} corresponds to the unit element, while the generators of degrees $1,2,3,$ and $4$ are
\begin{align*}
	&\sum_{n\geq 1}a_n\left(\frac{p_3}{2 n^6}-\frac{\left(H_n-H_{2 n}\right) p_0}{n^5}\right) = -\frac{\pi ^4}{60} \\
	& \sum_{n\geq 1}a_n \left( \frac{p_1}{n^7}-\frac{\left(3 H_n^{(2)}-H_{2 n}^{(2)}\right) p_0}{n^5}\right) = -2 \zeta (5) \\
	&\sum_{n\geq 1}a_n\left(-\frac{p_2}{n^8}+\frac{(H_n^{(3)}+H_{2 n}^{(3)}) p_0}{n^5}\right) = -2 \zeta (3)^2 \\
	&\sum_{n\geq 1}a_n\left(-\frac{49-245 n+320 n^2}{n^9}+\frac{(3 H_n^{(4)}+H_{2 n}^{(4)}) p_0}{2 n^5}\right) = -\zeta (7) 
\end{align*}
respectively, where, for brevity, we write $p_0 = 205 n^2-160 n+32,
p_1 = 553 n^2-440 n+88,
p_2 = 215 n^2-180 n+36,
	p_3 = 123 n^2-128 n+32$. The remaining generator, of degree $5$, is the $\zeta(5,3)$ series in \eqref{aux_00}. Non-generator identities look more complicated. For example, the remaining relation in degree $2$ is
$$\sum_{n\geq 1}a_n \left(-\frac{4-28 n+43 n^2}{n^7}+\frac{\left(2 H_n^2-4 H_n H_{2 n}+2 (H_{2 n})^2+H_n^{(2)}\right) p_0}{4 n^5}-\frac{\left(H_n-H_{2 n}\right) p_3}{2 n^6}\right) = -\frac{\zeta (5)}{2}.$$

The numerical evidence supporting this interpretation is substantial, as illustrated by all ten examples examined in the manuscript. Throughout, we highlight the identities corresponding to generators in \textcolor{NavyBlue}{blue}. The $S_5$ symmetry of a concrete case related to Example \ref{Ex_10} and the $\zeta(5)$ identity in \eqref{aux_7} was recently described in \cite{shalev2026symmetric}.

\subsection{The computational aspect}
This article is largely computational, so we briefly recall the computational setup. The WZ relation~\eqref{WZ} implies the summation formula
\begin{equation}\label{aux_8}\sum_{k\geq 0} F(0,k) = \sum_{n\geq 0} G(n,0) \qquad \text{(up to simpler terms)}.\end{equation}
In our applications, $F(0,k)$ arises from a WZ-seed. Typically, $\sum_{k\geq 0} F(0,k)$ is straightforward to handle, whereas $\sum_{n\geq 0} G(n,0)$ is the less straightforward and more interesting part. When the WZ-seed involves auxiliary parameters $a,b,c,\ldots$, these parameters also appear in $G(n,0)$; extracting appropriate coefficients then allows one to predict \textit{a priori} which series admit closed-form evaluations. \par

Examples in \cite{au2025multiple, chu2020alternating} show that the required computations can be enormous. Implementing the process is also far from straightforward: it requires extensive polynomial arithmetic, and an inefficient implementation can increase the running time by several orders of magnitude, preventing the computations needed for a complete proof from being carried out. \par

For this purpose, the author has developed a Mathematica package that automates this laborious coefficient-extraction process and completes it within a reasonable amount of time. The package allows us to discover many such series identities and, in some cases, to prove them rigorously, as explained below. The package and its source code can be downloaded from either of the following locations:
\begin{tcolorbox}
	Mathematica package \textsf{WZSeedTools}
	\tcblower
	\begin{itemize}[leftmargin = *]
		\item \url{https://www.researchgate.net/publication/400024737} or
		\item \url{https://sites.google.com/view/kc-au/2602-08721}
	\end{itemize}
\end{tcolorbox}\par

Our methodology allows us to discover many identities. We first enumerate the identities obtained by extracting coefficients, many of which involve harmonic numbers, and then ask whether there exists a linear combination in which relatively few harmonic-number terms remain. This approach yields a large collection of identities, some of which are presented in this article. Many others are too long to display here; readers interested in seeing them are encouraged to download and experiment with the Mathematica package mentioned above. \par

The paper is organized as follows. Section~2 provides background material on WZ-seeds and the technique of coefficient extraction. In Section~3, we present six examples in which the summation formula~\eqref{aux_8} can be written explicitly, allowing coefficient extraction to yield exact and rigorous evaluations; the conjecture above is also illustrated in the first three examples. Section~4 contains four further examples in which the full formula is not written out explicitly; instead, we focus on numerical identities obtained by extracting coefficients from $G(n,0)$.

		\section*{Acknowledgment}
The author is grateful to Dr. Erik Panzer for his valuable correspondence regarding this pre-print. The author has received funding from the European Research Council (ERC) under the European Union’s Horizon 2020 research and innovation programme (grant agreement No. 101001179). 

\section*{Declaration on AI usage}
The authors used ChatGPT 5.6 Sol for language polishing and TeX formatting of various rather long equations \textit{only}. All theoretical ideas and choice of examples come solely from the author. The accompanying Mathematica code is also solely written by the human author. 

\section{Preliminaries}
\subsection{Wilf-Zeilberger seeds}
We recall the definition of Wilf-Zeilberger seeds \cite[Section~3]{au2025wilf}:
\begin{definition}
	Let $f(a_1,\cdots,a_m,k)$ be a hypergeometric term\footnote{A hypergeometric term $f(n)$ is a function such that $f(n+1)/f(n)$ is a rational function in $n$.} in $a_i$ and $k$. It is called a \textit{WZ-seed} in the variables $a_i$ if, for all $A_i \in \mathbb{Z}$ and $K\in \mathbb{Z}$,
	$$F(n,k) = f(a_1+A_1 n,a_2+A_2 n,\cdots,a_m+A_m n,k+Kn)$$
	has a hypergeometric WZ-mate $G(n,k)$, that is, $$F(n+1,k) - F(n,k) = G(n,k+1) - G(n,k).$$
\end{definition}

Given a hypergeometric term $F(n,k)$, Gosper's algorithm can determine whether a hypergeometric $G(n,k)$ exists and, if so, find it \cite{gosper1978decision, paule1995mathematica}. For WZ-pairs arising from WZ-seeds, a more efficient algorithm\footnote{An implementation can be found in the Mathematica package mentioned in the introduction.} is available for finding $G(n,k)$, particularly when $A_i$ and $K$ are large \cite[Theorem~3.3]{au2025wilf}.


We reproduce below the list of WZ-seeds recorded in \cite{au2025wilf}, together with related instances from the literature in which these seeds have been used implicitly by various authors. This compilation also clarifies the origin of most WZ-pairs appearing in existing work. Moreover, every entry in our list admits a $q$-analogue; see \cite{au2024wilfQ}. \par

Some of the seeds below, especially those descended from \textsf{Dougall7F6} (\textsf{Seed3}, \textsf{Seed7}, and \textsf{Seed9}), yield new WZ-pairs and are the main ingredients in our proofs of the conjectural $1/\pi^3$ and $1/\pi^4$ formulas.

\begin{enumerate}
		\item  $\textsf{Gauss2F1}(a,b,c,k)$
		$$ = \frac{\Gamma (c-a) \Gamma (a+k) \Gamma (c-b) \Gamma (b+k)}{\Gamma (a) \Gamma (b) \Gamma (k+1) \Gamma (c+k) \Gamma (-a-b+c)}.$$
		Origin: Gauss $_2F_1$ summation formula. 
		
		\item  $\textsf{Dixon3F2}(a,b,c,k)$
		$$ = \frac{\Gamma (a-b+1) \Gamma (a-c+1) \Gamma (2 a+k) \Gamma (b+k) \Gamma (c+k) \Gamma (2 a-b-c+1)}{\Gamma (a) \Gamma (b) \Gamma (c) \Gamma (k+1) \Gamma (a-b-c+1) \Gamma (2 a-b+k+1) \Gamma (2 a-c+k+1)}.$$
		Origin: Dixon's $_3F_2$ summation formula. 
		
		\item $\textsf{Watson3F2}(a,b,c,k)$
		$$ = \frac{2^{-2 a-2 b+2 c} \Gamma \left(-a+c+\frac{1}{2}\right) \Gamma (2 a+k) \Gamma \left(-b+c+\frac{1}{2}\right) \Gamma (2 b+k) \Gamma (c+k)}{\Gamma (a) \Gamma (b) \Gamma (k+1) \Gamma (2 c+k) \Gamma \left(-a-b+c+\frac{1}{2}\right) \Gamma \left(a+b+k+\frac{1}{2}\right)}.$$ 
		Origin: Watson's $_3F_2$ summation formula. 
		
		\item $\textsf{Dougall5F4}(a,b,c,d,k)$
		$$ = \frac{(a+2 k) \Gamma (a+k) \Gamma (b+k) \Gamma (c+k) \Gamma (d+k) \Gamma (a-b-c+1) \Gamma (a-b-d+1) \Gamma (a-c-d+1)}{\Gamma (b) \Gamma (c) \Gamma (d) \Gamma (k+1) \Gamma (a-b+k+1) \Gamma (a-c+k+1) \Gamma (a-d+k+1) \Gamma (a-b-c-d+1)}.$$ 
		Origin: the very-well-poised $_5F_4$ summation formula. 
		
		\item $\textsf{Dougall7F6}(a,b,c,d,e,k)$
		$$ = \frac{\splitfrac{(-1)^{a+e} (a+2 k) \Gamma (a+k) \Gamma (b+k) \Gamma (c+k) \Gamma (d+k) \Gamma (e+k) \Gamma (a-b-c+1) \Gamma (a-b-d+1)}{ \Gamma (a-c-d+1) \Gamma (-a+b+c+e) \Gamma (-a+b+d+e) \Gamma (-a+c+d+e) \Gamma (2 a-b-c-d-e+k+1)}}{\splitfrac{\Gamma (b) \Gamma (c) \Gamma (d) \Gamma (e) \Gamma (k+1) \Gamma (-a+b+e) \Gamma (a-b+k+1) \Gamma (-a+c+e) \Gamma (a-c+k+1) \Gamma (-a+d+e) \Gamma (a-d+k+1)}{ \Gamma (a-e+k+1) \Gamma (a-b-c-d+1) \Gamma (2 a-b-c-d-e+1) \Gamma (-a+b+c+d+e+k)}}.$$ 
		Origin: the very-well-poised, $2$-balanced, terminating $_7F_6$ summation formula. 
		
		\item $\textsf{Balanced3F2}(a,b,c,d,k)$
		$$ = \frac{\Gamma (d-a) \Gamma (a+k) \Gamma (d-b) \Gamma (b+k) \Gamma (d-c) \Gamma (c+k) (-1)^{a+b+c-d}}{\Gamma (a) \Gamma (b) \Gamma (c) \Gamma (k+1) \Gamma (d+k) \Gamma (-a-b+d) \Gamma (-a-c+d) \Gamma (-b-c+d) \Gamma (a+b+c-d+k+1)}.$$ 
		Origin: $1$-balanced terminating $_3F_2$ summation formula.
		
		\item $\textsf{Dixon6F5}(a,b,c,k)$
		$$ = \frac{\splitfrac{2^{4 a-2 b-2 c}(-1)^{b+c+k} \Gamma(a+\frac{1}{2})  \Gamma (a-b+1) \Gamma (a-c+1)  \Gamma (-b+k+1) \Gamma (b+k) \Gamma (-c+k+1) }{\Gamma (c+k) (2 a-b-c+2 k+1) \Gamma (a-b-c+\frac{3}{2}) \Gamma (2 a-b-c+1) \Gamma (2 a-b-c+k+1)}}{\Gamma (k+1) \Gamma (2 a-b+k+1) \Gamma (2 a-c+k+1) \Gamma (2 a-b-2 c+k+2) \Gamma (2 a-2 b-c+k+2)}.$$ 
		Origin: combining Dixon's $_3F_2$ evaluation with Whipple's formula transforming an arbitrary $_3F_2$ at $1$ to a very-well-poised $_6F_5$ at $-1$.
		
		\item $\textsf{Seed1}(a,c,k)$
		$$ = \frac{2^{a-2 c-2 k} \Gamma (-a+2 c-1) \Gamma (a+2 k)}{\Gamma (a) \Gamma (k+1) \Gamma \left(-a+c-\frac{1}{2}\right) \Gamma (c+k)}.$$ 
		Origin: substitute $a \to a/2$ and $b\to a/2+1/2$ in \textsf{Gauss2F1}. 
		
		\item $\textsf{Seed2}(a,c,d,k)$
		$$ = \frac{2^{-2 c-2 k} (-1)^{a+c-d} \Gamma (-a+2 d-1) \Gamma (a+2 k) \Gamma (d-c) \Gamma (c+k)}{\Gamma (a) \Gamma (c) \Gamma (k+1) \Gamma \left(-a+d-\frac{1}{2}\right) \Gamma (d+k) \Gamma (-a-2 c+2 d-1) \Gamma \left(a+c-d+k+\frac{3}{2}\right)}.$$ 
		Origin: substitute $a \to a/2$ and $b\to a/2+1/2$ in \textsf{Balanced3F2}. 
		
		\item $\textsf{Seed3}(a,b,d,k)$
		$$ = \frac{2^{2 d} (a+2 k) \Gamma \left(a-b+\frac{1}{2}\right) \Gamma (a+k) \Gamma (b+2 k) \Gamma (d+k) \Gamma (2 a-b-2 d+1)}{\Gamma (b) \Gamma (d) \Gamma (k+1) \Gamma \left(a-b-d+\frac{1}{2}\right) \Gamma (2 a-b+2 k+1) \Gamma (a-d+k+1)}.$$ 
		Origin: substitute $b \to b/2$ and $c\to b/2+1/2$ in \textsf{Dougall5F4}. 
		
		\item $\textsf{Seed4}(a,b,k)$
		$$ = \frac{3^{b+1} (a+2 k) \Gamma (3 a-2 b) \Gamma (a+k) \Gamma (b+3 k)}{\Gamma (b) \Gamma (k+1) \Gamma (a-b) \Gamma (3 a-b+3 k+1)}.$$ 
		Origin: substitute $b \to b/3$, $c\to b/3+1/3$, and $d\to b/3+2/3$ in \textsf{Dougall5F4}. 
		
		\item $\textsf{Seed7}(a,b,d,e,k)$
		$$ = \frac{\splitfrac{(-1)^{a+e} (a+2 k) \Gamma \left(a-b+\frac{1}{2}\right) \Gamma (a+k) \Gamma (b+2 k) \Gamma (d+k) \Gamma (e+k) \Gamma (2 a-b-2 d+1)}{\Gamma \left(-a+b+e+\frac{1}{2}\right) \Gamma (-2 a+b+2 d+2 e) \Gamma \left(2 a-b-d-e+k+\frac{1}{2}\right)}}{\splitfrac{\Gamma (b) \Gamma (d) \Gamma (e) \Gamma (k+1) \Gamma \left(a-b-d+\frac{1}{2}\right) \Gamma (-2 a+b+2 e) \Gamma (2 a-b+2 k+1) \Gamma (-a+d+e) \Gamma (a-d+k+1)}{\Gamma (a-e+k+1) \Gamma \left(2 a-b-d-e+\frac{1}{2}\right) \Gamma \left(-a+b+d+e+k+\frac{1}{2}\right)}}.$$ 
		Origin: substitute $b \to b/2$ and $c\to b/2+1/2$ in \textsf{Dougall7F6}. 
		
		\item $\textsf{Seed9}(a,b,d,k)$
		$$ = \frac{\splitfrac{(a+2 k) (-1)^d 2^{4 a-2 b-2 d} \Gamma \left(a-b+\frac{1}{2}\right) \Gamma (a+k) \Gamma (b+2 k) \Gamma (d+2 k)}{ \Gamma (-2 a+2 b+d+1) \Gamma (-2 a+b+2 d+1) \Gamma (2 a-b-d+k)}}{\Gamma (b) \Gamma (d) \Gamma (k+1) \Gamma \left(-a+d+\frac{1}{2}\right) \Gamma (-2 a+b+d) \Gamma (2 a-b+2 k+1) \Gamma (2 a-d+2 k+1) \Gamma (-a+b+d+k+1)}.$$ 
		Origin: substitute $b \to b/2$, $c\to b/2+1/2$, $d\to d/2$, and $e\to d/2+1/2$ in \textsf{Dougall7F6}. 
		
		\item $\textsf{Seed10}(a,d,k)$
		$$ = \frac{(-1)^{a-d} 3^{-a-3 k} \Gamma (-a+3 d-2) \Gamma (a+3 k)}{\Gamma (a) \Gamma (k+1) \Gamma (-2 a+3 d-3) \Gamma (d+k) \Gamma (a-d+k+2)}.$$ 
		Origin: substitute $a\to a/3$, $b\to a/3+1/3$, and $c\to a/3+2/3$ in \textsf{Balanced3F2}. 
	\end{enumerate}
	
	Given a WZ-seed, one can construct WZ-pairs. The following proposition links such pairs to infinite series.
	\begin{proposition}\cite{au2025wilf} \label{WZ_prop}
		Suppose $F,G: \mathbb{N}^2\to \mathbb{C}$ are two functions satisfying
		$$F(n+1,k) - F(n,k) = G(n,k+1)-G(n,k)$$
		and the following conditions:
		\begin{itemize}
			\item $\sum_{k\geq 0} F(0,k)$ converges;
			\item $\sum_{n\geq 0} G(n,0)$ converges;
			\item $\lim_{k\to\infty} G(n,k) = g(n)$ exists for each $n\in \mathbb{N}$, and $\sum_{n\geq 0} g(n)$ converges.
		\end{itemize}
		Then $\lim_{n\to \infty} \sum_{k\geq 0} F(n,k)$ exists and is finite, and
		$$\sum_{k\geq 0} F(0,k) + \sum_{n\geq 0} g(n) = \sum_{n\geq 0} G(n,0) + \lim_{n\to \infty} \sum_{k\geq 0} F(n,k)$$
	\end{proposition}
	The simplest case is that in which both $g(n)$ and $\lim_{n\to \infty} \sum_{k\geq 0} F(n,k)$ vanish, yielding the formula
	$$\sum_{k\geq 0} F(0,k) = \sum_{n\geq 0} G(0,n).$$
	\begin{remark}
	There are many implementations on finding WZ-mate \cite{Koutschan09, paule1997mathematica, paule1995mathematica}. For WZ-seed, however, there is an even more efficient algorithm than the generic algorithm \cite[Theorem~3.3]{au2025wilf}. This is implemented in the Mathematica package mentioned in the introduction. 
	\end{remark}

\subsection{Harmonic numbers and coefficient extraction}
	We introduce notation that will be used throughout the paper. Let $f(a_1,\cdots,a_m)$ be a function analytic around the point $(a_1,\cdots,a_m) = (0,\cdots,0)$. We denote by $$[f(a_1,\cdots,a_m)][a_1^{r_1}\cdots a_m^{r_m}], \qquad r_i\in \mathbb{Z}_{\geq 0}$$
	the $a_1^{n_1}\cdots a_m^{n_m}$ coefficient in its power-series expansion. Let $$V_N(f) := \text{Span}_\mathbb{Q}\{[f(a_1,\cdots,a_m)][a_1^{r_1}\cdots a_m^{r_m}] | r_1 + \cdots + r_m = N\}$$
	be the $\mathbb{Q}$-span of the coefficients of all degree-$N$ monomials. If $f$ depends on an additional parameter $n\in \mathbb{N}$, then each $[f(a_1,\cdots,a_m)][a_1^{r_1}\cdots a_m^{r_m}]$ can be viewed as a function $\mathbb{N} \to \mathbb{C}$, and hence $V_N(f)$ can be viewed as a subspace of the $\mathbb{Q}$-vector space of functions from $\mathbb{N}\to \mathbb{C}$. The principal cases of interest are those in which $f(a_1,\cdots,a_m)$ is a product of rising factorials (i.e., Pochhammer symbols). \par
	
		For a fixed positive integer $N$ and $0<\gamma \leq 1$, define the harmonic sums
	$$H_N^{(s)}(\{\gamma\}^k) := \sum_{N-1 \geq n_1 >  \cdots > n_k \geq 0} \frac{1}{(n_ 1+\gamma)^{s}\cdots (n_k+\gamma)^{s}}, \qquad H_N^{(s)}(\gamma) := \sum_{n=0}^{N-1} \frac{1}{(n+\gamma)^s}.$$
	We recall the following facts (see, for example, \cite{au2025multiple} for more details).
	\begin{itemize}[leftmargin=*]
	\item When $k=\gamma=1$, $H_N^{(s)}(\{\gamma\}^k)$ reduces to the familiar harmonic number $H_N^{(s)} := 1+2^{-s} +\cdots + N^{-s}$. 
	\item When $k>1$, $H_N^{(s)}(\{\gamma\}^k)$ can be expressed as a polynomial in $H_N^{(i)}(\gamma), 1\leq i\leq s$, using Bell polynomials.
	\item When $\gamma = 1/2$, $H_N^{(s)}(1/2)$ can also be expressed in terms of familiar harmonic numbers as $2^s H_{2 N}^{(s)}-H_N^{(s)}$.
	\end{itemize}
We also have an easy-to-prove series expansion about $a=0$:
	\begin{equation}\label{pochhammerexpansion}(\gamma+a)_n = (\gamma)_n \left(1+ \sum_{k\geq 1} H_N^{(1)}(\{\gamma\}^k) a^k \right),\end{equation}
	which is the source of the harmonic numbers that arise when one extracts coefficients from a rising factorial. 
	\section{Explicit Examples}
	For Examples \ref{Ex_1} through \ref{Ex_4}, Proposition \ref{WZ_prop} yields the relation $$\sum_{n\geq 0} G(n,0) = \sum_{k\geq 0} F(0,k),$$
	since $g(n) = 0$ and $\lim_{n\to \infty} \sum_{k\geq 0} F(n,k) = 0$. All numerical identities listed in this section have been rigorously proved. We refer interested readers to the Mathematica package and files mentioned in the introduction for computational details. We reproduce the address here.
\begin{tcolorbox}
	Mathematica package \textsf{WZSeedTools}
	\tcblower
	\begin{itemize}[leftmargin = *]
		\item \url{https://www.researchgate.net/publication/400024737} or
		\item \url{https://sites.google.com/view/kc-au/2602-08721}
	\end{itemize}
\end{tcolorbox}
	\begin{romexample}\label{Ex_1}
	Set $$F(n,k) = \textsf{Dougall5F4}\left(a, b, c, 1 + d - n, k +e + n\right).$$ Then $$\sum_{n\geq 1} G(n-1,0) = \sum_{k\geq 0} F(0,k)$$ yields the identity \begin{multline}\label{aux_9} \sum_{n\geq 1} \frac{(-1)^{n-1} (1-d)_{n-1} (a+e+1)_{n-1} (b+e+1)_{n-1} (c+e+1)_{n-1} (a-b-d+1)_{n-1} (a-c-d+1)_{n-1} P(n)}{(e+1)_n (a-b+e+1)_n (a-c+e+1)_n (a-d+e+1)_{2 n} (a-b-c-d+1)_n}  \\ = \sum_{k\geq 0} \frac{(a+2 e+2 k+2) (a+e+1)_k (b+e+1)_k (c+e+1)_k (d+e+1)_k}{(e+1)_{k+1}(a-b+e+1)_{k+1} (a-c+e+1)_{k+1} (a-d+e+1)_{k+1}} \end{multline}
	where $P(n)$ is the polynomial \begin{multline*}P(n) = a^3-a^2 (b+c+2 d-3 e)+a \left(b (c+d-2 e)+c d-2 c e+d^2-5 d e+3 e^2\right) \\ +n \left(5 a^2-a (3 b+3 c+7 d-11 e)+b (c+2 d-4 e)+2 c d-4 c e+2 d^2-8 d e+5 e^2\right) \\ +n^2 (9 a-3 (b+c+2 d-3 e))+e \left(b (c+2 d-e)+2 c d-c e+2 d^2-2 d e+e^2\right)+5 n^3.\end{multline*}
If we define $$a_n = \left(\frac{-1}{4}\right)^n \frac{(1)_n}{(\frac12)_n}$$ and denote the summand in $n$ by $S := S(a,b,c,d,e;n)$, then $$[S][a^0b^0c^0d^0e^0] = a_n \frac{5}{n^3} \implies \dim V_0(S) = 1.$$
There are five monomials of degree $1$; their respective coefficients in $S$ are:
	$$\begin{aligned}
	[S][a^1b^0c^0d^0e^0] &= a_n\left(-\frac{5 H_{2 n}}{n^3}-\frac{6}{n^4}\right), \\
	[S][a^0b^1c^0d^0e^0] &= a_n\left(\frac{10 H_n}{n^3}-\frac{3}{n^4}\right), \\
	[S][a^0b^0c^1d^0e^0] &= a_n\left(\frac{10 H_n}{n^3}-\frac{3}{n^4}\right), \\
	[S][a^0b^0c^0d^1e^0] &= a_n\left(-\frac{10 H_n}{n^3}+\frac{5 H_{2 n}}{n^3}+\frac{9}{n^4}\right), \\
	[S][a^0b^0c^0d^0e^1] &= a_n\left(-\frac{5 H_{2 n}}{n^3}-\frac{6}{n^4}\right).
	\end{aligned}$$
However, they span only a two-dimensional space: the last three coefficients are linear combinations of the first two, so $\dim V_1(S) = 2$. Similarly, the $15$ coefficients of degree-$2$ monomials span only a four-dimensional space, generated by:
\begin{multline*}a_n\left(\frac{H_n^{(2)}}{n^3}-\frac{9}{5 n^5}\right),\qquad a_n\left(-\frac{3 H_n}{5 n^4}+\frac{(H_n)^2}{n^3}+\frac{1}{2 n^5}\right), \\ a_n\left(\frac{6 H_n}{5 n^4}-\frac{3 H_{2 n}}{10 n^4}+\frac{H_{2 n} H_n}{n^3}+\frac{4}{5 n^5}\right),\qquad  a_n\left(\frac{12 H_{2 n}}{5 n^4}+\frac{(H_{2 n})^2}{n^3}+\frac{H_{2 n}^{(2)}}{n^3}+\frac{16}{5 n^5}\right).\end{multline*}
The table lists the dimension of $V_N(S)$ for the first few values of $N$. 
			\begin{table}[h]
			\begin{tabular}{|c|c|c|c|c|c|c|c|c|}
				\hline
				$N$                          & $0$ & $1$ & $2$  & $3$  & $4$  & $5$  & $6$ & $7$   \\ \hline
				$\dim V_N(S)$              & $1$ & $2$ & $4$ & $7$ & $12$ & $18$ & $27$ & $38$ \\ \hline
			\end{tabular}
			\caption{\small Number of independent summation formulas that can be generated by $S(a,b,c,d,e;n)$.}
		\end{table}

\begin{conjecture}
The generating function
$$\sum_{N\geq 0} \dim V_N(S) t^N = 1+2t+4t^2+7t^3+12t^4+18t^5+\cdots \stackrel{?}{=} \frac{1}{(1-t)^2 \left(1-t^2\right) \left(1-t^3\right) \left(1-t^4\right)}.$$
\end{conjecture}

Expanding the right-hand side of the formula above as a power series in $a,b,c,d,e$, we find that the coefficients of degree $N$ are MZVs of weight $N+3$. We therefore have a $\mathbb{Q}$-linear map: $$\Sigma: V_N(S) \to \MZV{}{N+3},\qquad f(n) \mapsto \sum_{n\geq 1} f(n).$$
The image of $V_0(S)$ is the well-known Apéry sum: $$\sum_{n\geq 0} a_n \frac{5}{n^3} = -2\zeta(3).$$
The image of $V_1(S)$ is
\begin{equation}\label{aux_1}\sum_{n\geq 1} a_n \left(\frac{H_n}{n^3}-\frac{3}{10 n^4}\right) = -\frac{\pi^4}{300},\qquad \sum_{n\geq 1} a_n \left(\frac{H_{2 n}}{n^3}+\frac{6}{5 n^4}\right) = -\frac{\pi^4}{75}.\end{equation}
The image of $V_2(S)$ is
\begin{align}\label{aux_2}
	\begin{split}&\sum_{n\geq 1}a_n\left(\frac{H_n^{(2)}}{n^3}-\frac{9}{5 n^5}\right) = \frac{2}{5}\zeta(5) \\
&\sum_{n\geq 1}a_n\left(-\frac{3 H_n}{5 n^4}+\frac{(H_n)^2}{n^3}+\frac{1}{2 n^5}\right) = \frac{-2}{5}\zeta(5)\\
&\sum_{n\geq 1}a_n\left( \frac{6 H_n}{5 n^4}-\frac{3 H_{2 n}}{10 n^4}+\frac{H_{2 n} H_n}{n^3}+\frac{4}{5 n^5}\right) = -\frac{7}{5}\zeta(5) \\
&\sum_{n\geq 1}a_n\left(\frac{12 H_{2 n}}{5 n^4}+\frac{(H_{2 n})^2}{n^3}+\frac{H_{2 n}^{(2)}}{n^3}+\frac{16}{5 n^5}\right) = \frac{-24}{5}\zeta(5)
\end{split}\end{align}
Therefore, we see\footnote{Assuming the standard conjecture on the $\mathbb{Q}$-dimension of $\MZV{}{N}$.} that the images of $V_0(S), V_1(S), V_2(S)$ under $\Sigma$ are all one-dimensional, spanned by $\zeta(3)$, $\zeta(4)$, and $\zeta(5)$, respectively. Several such series were conjectured by Sun \cite{sun2021book}, and some were proved in \cite{chu2020alternating} using the same formula and method as in this example. We refer the reader to \cite{chu2020alternating} for more explicit numerical examples. 
\begin{table}[h]
	\begin{tabular}{|c|c|}
		\hline
		$N$     & Spanning set of $\Sigma(V_N(S))$                                          \\ \hline
		$0,1,2$ & $\{\zeta(3+N)\}$                                                          \\ \hline
		$3$     & $\{\zeta(3)^2,\zeta(6)\}$                                                 \\ \hline
		$4$     & $\{\zeta(7),\zeta(3)\zeta(4)\}$                                           \\ \hline
		$5$     & $\{\zeta(5,3),\zeta(3)\zeta(5),\zeta(8)\}$                                \\ \hline
		$6$     & $\{\zeta (9),\zeta(4)\zeta (5),\zeta (3)^3,\zeta(6)\zeta (3)\}$           \\ \hline
		$7$     & $\{\zeta(7,3),\zeta(3)\zeta(7),\zeta(5)^2,\zeta(4)\zeta(3)^2,\zeta(10)\}$ \\ \hline
	\end{tabular}
\end{table}

\begin{conjecture}
	The generating function $$\sum_{N\geq 0} \dim \Sigma(V_N(S)) t^N \stackrel{?}{=} 1+t+t^2+2t^3+2t^4+\cdots$$ is a rational function in $t$.
\end{conjecture}
A step-by-step computational guide for this example is available in the Mathematica package mentioned in the introduction. 
	\end{romexample}
	
\begin{remark}
This simple example also has an interpretation in terms of colored MZVs, in which multiple Landen values arise naturally; see \cite{borwein2001central, broadhurst2015multiple, au2022iterated, xu2022sun}. However, the other examples in this article do not offer such an alternative interpretation, nor does this interpretation make the evaluations any easier. 
\end{remark}

\begin{romexample}\label{Ex_2}
Set $$F(n,k) = \textsf{Dougall5F4}\left(a,b-n+1,c,d+\frac{1}{2},k+n\right).$$ Then $$\sum_{n\geq 1} G(n-1,0) = \sum_{k\geq 0} F(0,k)$$ yields the identity 
		$$\begin{aligned} \sum_{n\geq 1} \frac{(-1)^{n-1} P(n) \times (1-b)_{n-1} (a+e+1)_{n-1} (c+e+1)_{n-1} (a-b-c+1)_{n-1} (d+e+\frac{1}{2})_n  (a-b-d+\frac{1}{2})_{n-1}}{2(e+1)_n (a-c+e+1)_n (a-d+e+\frac{1}{2})_n (a-b-c-d+\frac{1}{2})_n  (a-b+e+1)_{2 n}} & \\ = \sum_{k\geq 0} \frac{(a+2 e+2 k+2) (a+e+1)_k (b+e+1)_k (c+e+1)_k (d+e+\frac{1}{2})_{k+1}}{(e+1)_{k+1} (a-b+e+1)_{k+1} (a-c+e+1)_{k+1} (a-d+e+\frac{1}{2})_{k+1}}. \end{aligned},$$
where $P(n)$ is the following polynomial:
	\begin{multline*}
	P(n) = 2 a^3-a^2 (4 b+2 c+2 d-6 e+1)+a \left(2 b^2+b (2 c+2 d-10 e+1)+2 c d-4 c e+c-4 d e+6 e^2-2 e\right) \\ +e \left(4 b^2+b (4 c+4 d-4 e+2)+(c-e) (2 d-2 e+1)\right) \\ +n \left(10 a^2-a (14 b+6 c+6 d-22 e+3)+4 b^2+2 b (2 c+2 d-8 e+1)+2 c d-8 c e+c-8 d e+10 e^2-4 e\right) \\ + 3 n^2 (6 a-4 b-2 c-2 d+6 e-1)+10 n^3.
	\end{multline*}
Let $S_1 := S_1(a,b,c,d,e;n)$ denote the summand in $n$, let $S_2 := a S(a,b,c,d,e;n)$, where $S$ is defined as in the previous example, and set $$a_n := \left(\frac{-1}{4}\right)^n \frac{(1)_n}{(\frac12)_n}.$$ It can be shown, for example by using a transformation formula for $_7F_6$ \cite{wei2023two}, that each power-series coefficient of a degree-$N$ monomial on the right-hand side lies in $\MZV{}{N+2}$. For the map $\Sigma$ defined in the previous example by $f(n) \mapsto \sum_{n\geq 1}f(n)$, we then have
$$\Sigma: V_N(S_i) \to \MZV{}{N+2},\qquad \text{ for } i=1,2.$$
The dimension of $V_0(S_1)$ is one; its image under $\Sigma$ gives the identity
$$\sum_{n\geq 1} a_n \frac{10 n-3}{n^2 (2 n-1)} =  -\frac{\pi^2}{3}.$$
The dimension of $V_1(S_1)$ is three, and its image gives the identities
\begin{align*}
&\sum_{n\geq 1} a_n \frac{5}{n^3} = -2\zeta(3), \\
&\sum_{n\geq 1} a_n \left( \frac{(10 n-3) H_n}{3 n^2 (2 n-1)}-\frac{4 (3 n-1)}{3 n^2 (2 n-1)^2}\right) =  \frac{\zeta(3)}{5}, \\
&\sum_{n\geq 1} a_n \left( \frac{(10 n-3) H_{2 n}}{3 n^2 (2 n-1)}-\frac{2}{3 n^2 (2 n-1)}\right) =  -\frac{16\zeta(3)}{15}
\end{align*}
The first identity has already appeared in the previous example. We list the dimensions of $V_N(S_1)$, $V_N(S_2)$, and their sum in the following table. The table shows that, when $N\leq 4$, $V_N(S_2) \subset V_N(S_1)$, but $S_2$ gives additional relations when $N\geq 5$.

	\begin{table}[h]
	\begin{tabular}{|c|c|c|c|c|c|c|c|}
		\hline
		$N$                          & $0$ & $1$ & $2$  & $3$  & $4$  & $5$  & $6$   \\ \hline
		$\dim V_N(S_1)$              & $1$ & $3$ & $7$ & $15$ & $29$ & $49$ & $79$ \\ \hline
		$\dim V_N(S_2)$              & $0$ & $1$ & $2$  & $4$ & $7$ & $12$ & $18$  \\ \hline
		$\dim (V_N(S_1) + V_N(S_2))$ & $1$ & $3$ & $7$ & $15$ & $29$ & $52$ & $86$ \\ \hline
	\end{tabular}
	\caption{\small Number of independent identities that can be generated by $S_1(a,b,c,d;n)$ and $S_2(a,b,c,d;n)$, separately and in combination. The second row is identical to $\dim V_{N-1}(S)$ in the previous example.}
\end{table}

The image $\Sigma(V_N(S_1))$ coincides with $\MZV{}{N+2}$ when $N\leq 6$. However, this cannot hold for all $N$, since $\dim V_N(S_1) = O(N^4)$ while $\dim \MZV{}{N+2}$ is conjectured to grow exponentially. Nonetheless, we make the following conjecture.

\begin{conjecture}
	(a) The generating function $$\sum_{N\geq 0} \dim V_N(S_1) t^N \stackrel{?}{=} \frac{1}{(1-t)^3 (1-t^2)^2}-\frac{t^2}{1-t}.$$ 
	(b) The generating function $\sum_{N\geq 0} \dim  \Sigma(V_N(S_1)) t^N$ is a rational function in $t$. 
\end{conjecture}
Below, we list some identities arising from $V_N(S_1) + V_N(S_2)$, restricting attention to the simplest-looking ones:
\begin{align*}
	&\sum_{n\geq 1} a_n \left( \frac{(10 n-3) H_n^{(2)}}{3 n^2 (2 n-1)}-\frac{18 n-7}{3 n^4 (2 n-1)}\right) =  -\frac{41\pi^4}{2700}\\
	&\sum_{n\geq 1} a_n \left( \frac{(10 n-3) H_{2 n}^{(2)}}{3 n^2 (2 n-1)}-\frac{9 n-4}{3 n^4 (2 n-1)} \right) = -\frac{19\pi^4}{900}\\
	&\sum_{n\geq 1} a_n \left( \frac{(10 n-3) (H_n)^2}{3 n^2 (2 n-1)}-\frac{8 (3 n-1) H_n}{3 n^2 (2 n-1)^2}+\frac{4 \left(34 n^3-6 n^2-7 n+2\right)}{15 n^4 (2 n-1)^3} \right) =  -\frac{\pi^4}{300}\\
	&\sum_{n\geq 1} a_n \left( \frac{(10 n-3) \left(H_{2 n}\right){}^2}{3 n^2 (2 n-1)}-\frac{4 H_{2 n}}{3 n^2 (2 n-1)}+\frac{37 n-16}{15 n^4 (2 n-1)} \right) = -\frac{19\pi^4}{400}\\
	&\sum_{n\geq 1} a_n \left( \frac{(10 n-3) H_n^{(3)}}{3 n^2 (2 n-1)}-\frac{38 n-17}{3 n^5 (2 n-1)} \right) = -\frac{16}{3}\zeta(5) + \frac{2}{3}\pi^2\zeta(3)\\
	&\sum_{n\geq 1} a_n \left( \frac{(10 n-3) H_{2 n}^{(3)}}{3 n^2 (2 n-1)}+\frac{32 n-17}{12 n^5 (2 n-1)} \right) = \frac{8}{3}\zeta(5)  - \frac{7}{18}\pi^2\zeta(3)\\
	&\sum_{n\geq 1} a_n\left( \frac{H_{2 n}^{(2)}}{n^3}-\frac{32}{5 n^5}\right) = -\frac{41}{5}\zeta(5)+\frac{14}{15}\pi^2\zeta(3) \\
	&\sum_{n\geq 1} a_n \left(\frac{(10 n-3) H_n^{(4)}}{3 n^2 (2 n-1)}-\frac{4 (7 n-3)}{3 n^6 (2 n-1)} \right) = \frac{\pi^6}{1215} + \frac{8}{15}\zeta(3)^2\\
	&\sum_{n\geq 1} a_n \left(\frac{H_{2 n}^{(3)}}{n^3}+\frac{4 (10 n-3) H_{2 n}^{(4)}}{5 n^2 (2 n-1)}+\frac{6 n-5}{20 n^6 (2 n-1)} \right) = -\frac{3 \zeta (3)^2}{10}-\frac{23 \pi ^6}{7560}\\
	&\sum_{n\geq 1} a_n \left(\frac{10 H_n H_{2 n}^{(2)}}{n^3}-\frac{64 H_n}{n^5}-\frac{3 H_{2 n}^{(2)}}{n^4}-\frac{10 H_{2 n}^{(3)}}{n^3}+\frac{2}{n^6} \right) = \frac{23 \pi ^6}{378}-\frac{92 \zeta (3)^2}{5} \\
	&\sum_{n\geq 1} a_n\Bigg(-\frac{2 (-15+34 n)}{n^8 (-1+2 n)}+\frac{2 (-11+26 n) H_n^{(2)}}{n^6 (-1+2 n)}-\frac{(-7+18 n) ((H_n^{(2)})^2-H_n^{(4)})}{n^4 (-1+2 n)} \\ &\quad +\frac{(-3+10 n) ((H_n^{(2)})^3-3 H_n^{(2)} H_n^{(4)}+2 H_n^{(6)})}{3 n^2 (-1+2 n)} \Bigg) = \frac{2 \pi ^8}{4725}
\end{align*}
A step-by-step computational guide for this example is available in the Mathematica package mentioned in the introduction. Readers who wish to see more formulas can consult the package. \par

Although the summands in the seventh and tenth identities contain no factor $2n-1$ in their denominators, they cannot be derived from the previous example alone. This illustrates the importance of having more than one formula (in this case, $S_1$ and $S_2$) available when deriving identities involving harmonic numbers. 
\end{romexample}

\begin{romexample}\label{Ex_3}
	This example concerns a harmonic extension of the series \cite{amdeberhan1998hypergeometric}
	$$\sum_{n\geq 1} \left(\frac{-1}{2^{10}}\right)^n \frac{(1)_n^5}{(\frac12)_n^5}  \frac{32-160 n+205 n^2}{n^5} = -2\zeta(3).$$
	It is known that this series has a five-parameter extension \cite[Example~VI]{au2025multiple}; see also \cite{chu2014accelerating}:
\begin{multline*}\sum_{n\geq 1}(-1)^{n-1} P(n) \frac{\splitfrac{(a+1)_{n-1} (b+1)_{n-1}(c+1)_{n-1} (d+1)_{n-1} (-b+e+1)_{n-1} (-c+e+1)_{n-1} }{(-d+e+1)_{n-1} (a-b-c+e+1)_{n-1}(a-b-d+e+1)_{n-1} (a-c-d+e+1)_{n-1}}}{(e+1)_{2 n}(a-b+e+1)_{2 n}(a-c+e+1)_{2 n} (a-d+e+1)_{2 n} (a-b-c-d+2 e+1)_{2 n}} \\ = \sum_{k\geq 0}\frac{(a+e+2 k+2) (a+1)_k (b+1)_k (c+1)_k (d+1)_k}{(e+1)_{k+1} (a-b+e+1)_{k+1} (a-c+e+1)_{k+1} (a-d+e+1)_{k+1}},
\end{multline*}
where $P(n)\in \mathbb{Z}[a,b,c,d,e,n]$ is a lengthy polynomial whose expression can be found in the Mathematica file. Let $S:= S(a,b,c,d,e;n)$ denote the summand in $n$, and set $a_n :=  \left(\frac{-1}{2^{10}}\right)^n \frac{(1)_n^5}{(\frac12)_n^5}$. Each power-series coefficient of a degree-$N$ monomial on the right-hand side lies in $\MZV{}{N+3}$. 
	\begin{table}[h]
		\begin{tabular}{|c|c|c|c|c|c|c|c|c|}
			\hline
			$N$                          & $0$ & $1$ & $2$  & $3$  & $4$  & $5$  & $6$ & $7$   \\ \hline
			$\dim V_N(S)$              & $1$ & $1$ & $2$ & $3$ & $5$ & $7$ & $10$ & $13$ \\ \hline
		\end{tabular}
		\caption{\small Number of independent summation formulas that can be generated by $S(a,b,c,d,e;n)$.}
		\label{tab:my-table}
	\end{table}
	\begin{conjecture}\label{gen_conj_Ex3}
		The generating function
		$$\sum_{N\geq 0} \dim V_N(S) t^N \stackrel{?}{=} \frac{1}{(1-t)(1-t^2)(1-t^3)(1-t^4)(1-t^5)}.$$
	\end{conjecture}
We list some examples of identities arising from $V_N(S)$; for $N\geq 3$, we include only the simplest-looking ones. To abbreviate the formulas, we write
	$$p_0 = 205 n^2-160 n+32, \quad
		p_1 = 553 n^2-440 n+88, \quad
		p_2 = 215 n^2-180 n+36, \quad
		p_3 = 123 n^2-128 n+32.$$
	\begin{align*}
		&\color{NavyBlue}\sum_{n\geq 1}a_n\left(\frac{p_3}{2 n^6}-\frac{\left(H_n-H_{2 n}\right) p_0}{n^5}\right) = -\frac{\pi ^4}{60} \\
		&\color{NavyBlue} \sum_{n\geq 1}a_n \left( \frac{p_1}{n^7}-\frac{\left(3 H_n^{(2)}-H_{2 n}^{(2)}\right) p_0}{n^5}\right) = -2 \zeta (5) \\
		&\sum_{n\geq 1}a_n \left(-\frac{4-28 n+43 n^2}{n^7}+\frac{\left(2 H_n^2-4 H_n H_{2 n}+2 (H_{2 n})^2+H_n^{(2)}\right) p_0}{4 n^5}-\frac{\left(H_n-H_{2 n}\right) p_3}{2 n^6}\right) = -\frac{\zeta (5)}{2} \\
		&\color{NavyBlue} \sum_{n\geq 1}a_n\left(-\frac{p_2}{n^8}+\frac{(H_n^{(3)}+H_{2 n}^{(3)}) p_0}{n^5}\right) = -2 \zeta (3)^2 \\
		&\color{NavyBlue}\sum_{n\geq 1}a_n\left(-\frac{49-245 n+320 n^2}{n^9}+\frac{(3 H_n^{(4)}+H_{2 n}^{(4)}) p_0}{2 n^5}\right) = -\zeta (7) \\
		&\color{NavyBlue} \sum_{n\geq 1} a_n \left(\frac{351-1755 n+2240 n^2}{n^{10}}-\frac{\left(11 H_n^{(5)}-H_{2 n}^{(5)}\right) p_0}{n^5}\right) = -4 \zeta(5,3)-\frac{\pi ^8}{4725} \\
		&\sum_{n\geq 1}a_n\left(\frac{68-340 n+421 n^2}{n^{11}}+\frac{((H_n^{(3)})^2-3 H_n^{(6)}+2 H_n^{(3)} H_{2 n}^{(3)}+(H_{2 n}^{(3)})^2+H_{2 n}^{(6)}) p_0}{2 n^5}-\frac{(H_n^{(3)}+H_{2 n}^{(3)}) p_2}{n^8}\right) \\ &\quad = -\frac{4\zeta(3)^3}{3}-\frac{2 \zeta (9)}{3}.
	\end{align*}
The simple pattern of identities involving only $H_n^{(N)}$ and $H_{2n}^{(N)}$ ends at $N=6$; at that point, the simplest form is the last equation. The first four colored equations were conjectured by Sun \cite{sun2022conjectures} and proved in \cite{wei2023some, au2025multiple}. \par 

A step-by-step computational guide for this example is available in the Mathematica package mentioned in the introduction. Readers who wish to see more formulas can consult the package.
\end{romexample}

\begin{remark}
As mentioned in the introduction, the generating function of $\dim V_N(S)$ suggests that $V_N(S)$ can be interpreted as the graded pieces of an algebra with five generators of degrees $1$, $2$, $3$, $4$, and $5$. It is natural to ask which terms play the role of these generators. A simple heuristic is that they are precisely the {\color{NavyBlue}“colored”} terms. In other words,
\begin{equation}\label{aux_10}
	\text{generators of this algebra} \Longleftrightarrow \text{summands without products of harmonic numbers}.
\end{equation}
Although this interpretation is somewhat speculative, we will see that this correspondence holds uniformly across a wide range of examples (see Examples~\ref{Ex_1}, \ref{Ex_5}, \ref{Ex_8}, \ref{Ex_9}, and~\ref{Ex_10}), where we shall also highlight these purported generators.
\end{remark}

\begin{romexample}\label{Ex_4}
	This example concerns a harmonic extension of the series \cite[Example~24]{chu2014accelerating}
	\begin{equation}\label{aux_3}\sum_{n\geq 1} \left(\frac{-1}{27}\right)^n \frac{(1)_n^2}{(\frac13)_n (\frac23)_n}  \frac{7 n-2}{n^2 (2 n-1)} = -\frac{\pi ^2}{12}.\end{equation}
Several WZ-seeds provide multivariable extensions of this series, namely,
$$\begin{aligned}
F_1(n,k) &= \textsf{Dougall5F4}\left(a+n+\frac{1}{2},b+\frac{1}{2},c,d-n,k+d+n+\frac{1}{2}\right)\\
F_2(n,k) &= \textsf{Dougall5F4}\left(a,b,c,d-n,e+k+2n+1\right)\\
F_3(n,k) &= \textsf{Seed3}\left(a,b-n,c,k+d+n+1\right)
\end{aligned}$$
The three summation formulas are then $$\sum_{k\geq 0} F_i(0,k) = \sum_{n\geq 1} G_i(n-1,0).$$ The first two are
\begin{multline*}\sum_{n\geq 1} \frac{\splitfrac{(-1)^{n-1} P_1(n) (1-d)_{n-1} (a+e+1)_{2 n-2} (b+e+1)_{n-1} (c+e+\frac{1}{2})_{n-1} (a-b-c+1)_{n-1}}{\times (a-b-d+1)_{2 n-2} (a-c-d+\frac{1}{2})_{2 n}}}{(e+\frac{1}{2})_n (2 a-2 c-2 d+4 n-1) \left(a-b+e+\frac{1}{2}\right)_{2 n} (a-c+e+1)_{2 n} (a-b-c-d+1)_{2 n} (a-d+e+1)_{3 n}}  \\=\sum_{k\geq 0}\frac{(2 a-2 c-2 d+1) (2 a+4 e+4 k+3) (a+e+1)_k (b+e+1)_k (c+e+\frac{1}{2})_k (d+e+\frac{1}{2})_k}{2 (e+\frac{1}{2})_{k+1} (a-b+e+\frac{1}{2})_{k+1} (a-c+e+1)_{k+1} (a-d+e+1)_{k+1}},
\end{multline*}
and
\begin{multline*}
\sum_{n\geq 1} \frac{\splitfrac{(-1)^{n-1} P_2(n) (1-d)_{n-1} (a+e+1)_{2 n-2} (b+e+1)_{2 n-2} (c+e+1)_{2 n-2}}{ (d+e+1)_{n-1} (a-b-d+1)_{n-1} (a-c-d+1)_{n-1}}}{(e+1)_{2 n} (a-b+e+1)_{2 n} (a-c+e+1)_{2 n} (a-d+e+1)_{3 n} (a-b-c-d+1)_n}  = \text{the RHS of }\eqref{aux_9}.\end{multline*}
Here $P_1(n), P_2(n) \in \mathbb{Z}[a,b,c,d,e,n]$ are lengthy polynomials\footnote{They can also be found in \cite{chu2014accelerating} in a shorter, unexpanded form.} whose expressions can be found in the data file. The third formula is shorter:
\begin{multline*}
	\sum_{n\geq 1}\frac{(-1)^{n-1} P_3(n) (1-b)_{n-1} \left(a-b+\frac{1}{2}\right)_{n-1} (a+d+1)_{n-1} (b+2 d+1)_n (c+d+1)_{n-1} (2 a-b-2 c+1)_{n-1}}{2 (d+1)_n \left(a-b-c+\frac{1}{2}\right)_n (2 a-b+2 d+1)_{3 n} (a-c+d+1)_n} \\
	=\sum_{k\geq 0} \frac{(a+2 d+2 k+2) (a+d+1)_k (b+2 d+1)_{2 k+1} (c+d+1)_k}{(d+1)_{k+1} (2 a-b+2 d+1)_{2 k+2} (a-c+d+1)_{k+1}}
\end{multline*}
where \begin{multline*}P_3(n) = 4 a^3 - 2 a^2 (1 + 3 b + 2 c - 6 d) + 
	2 d (b + 2 b^2 + c + 3 b c - d - 3 b d - 2 c d + 2 d^2) + 
	a (b + 2 b^2 + 2 c + 4 b c - 4 d - 14 b d - 8 c d + 12 d^2) \\ + [18 a^2 - a (5 + 18 b + 12 c - 38 d) + 2 (b + 2 b^2 + c + 3 b c - 3 d - 10 b d - 7 c d + 9 d^2)] n + 
	2 (-2 + 14 a - 7 b - 5 c + 14 d) n^2 + 14 n^3.\end{multline*}
When $a=b=c=d=e=0$, the first and third formulas specialize to \eqref{aux_3}, whereas the second specializes to the identity in \cite[Example~21]{chu2014accelerating}: $$\sum_{n\geq 1} a_n \frac{5-32 n+56 n^2}{n^3 (-1+2 n)^2} = -4 \zeta (3), \qquad a_n :=  \left(\frac{-1}{27}\right)^n \frac{(1)_n^2}{(\frac13)_n (\frac23)_n}.$$
The coefficients of degree-$N$ monomials in the first and third formulas take values in $\MZV{}{2+N}$, whereas those in the second formula take values in $\MZV{}{3+N}$. Let $S_1 := S_1(a,b,c,d,e;n)$ denote the summand in the first formula, let $S_2 := a\times S_2(a,b,c,d,e;n)$ denote the summand in the second formula, let $S_3 :=  S_3(a,b,c,d;n)$ denote the summand in the third formula. The coefficients of degree-$N$ monomials on the right-hand sides of all three formulas are elements of $\MZV{}{N+2}$. 
	\begin{table}[h]
	\begin{tabular}{|c|c|c|c|c|c|c|c|}
		\hline
		$N$                          & $0$ & $1$ & $2$  & $3$  & $4$  & $5$  & $6$   \\ \hline
		$\dim V_N(S_1)$              & $1$ & $4$ & $11$ & $24$ & $46$ & $80$ & $130$ \\ \hline
		$\dim V_N(S_2)$              & $0$ & $1$ & $3$ & $6$ & $13$ & $23$ & $38$ \\ \hline
		$\dim V_N(S_3)$              & $1$ & $3$ & $7$  & $13$ & $22$ & $34$ & $50$  \\ \hline
		$\dim (V_N(S_1)+V_N(S_2))$              & $1$ & $4$ & $12$  & $27$ & $54$ & $96$ & $159$  \\ \hline
		$\dim (V_N(S_1) + V_N(S_3))$ & $1$ & $4$ & $11$ & $26$ & $52$ & $93$ & $153$ \\ \hline
		$\dim (V_N(S_2)+V_N(S_3))$              & $1$ & $3$ & $8$  & $16$ & $30$ & $50$ & $79$  \\ \hline
		$\dim (V_N(S_1)+V_N(S_2)+V_N(S_3))$              & $1$ & $4$ & $12$  & $29$ & $59$ & $108$ & $180$  \\ \hline
	\end{tabular}
	\label{Table_Ex_4}
	\caption{\small Number of independent summation formulas that can be generated by $S_1$, $S_2$, and $S_3$, separately and in combination.}
\end{table}
	\begin{conjecture}
	The generating function for each row of the table above is a rational function. In particular,
	$$\sum_{N\geq 0} \dim V_N(S_1) t^N \stackrel{?}{=} \frac{1}{(1-t)^4 (1-t^2)}, \qquad \sum_{N\geq 0} \dim V_N(S_3)t^N  \stackrel{?}{=}  \frac{1}{(1-t)^3 (1-t^2)}.$$
\end{conjecture}

The table shows that, when $N\leq 2$, the first formula subsumes the second; for $N\geq 3$, however, the second and third formulas provide new relations. One such identity at $N=3$, which cannot be proved from $S_1$ alone, was conjectured by Sun\footnote{\url{https://mathoverflow.net/questions/486353}}:
$$\sum_{n\geq 1} a_n \frac{17-136 n+408 n^2-640 n^3+560 n^4}{n^5 (-1+2 n)^4} = 180 \zeta (5)-\frac{56 \pi ^2 \zeta (3)}{3}.$$
We list some examples of identities from $V_N(S_1)+V_N(S_2)+V_N(S_3)$, including only the simplest-looking ones. For explicit and rigorous computations of the right-hand sides, we refer readers to the accompanying Mathematica file\footnote{The file is available at \url{https://www.researchgate.net/publication/400024737}.}.
	\begingroup
	\let\originalaligned\aligned
	\let\endoriginalaligned\endaligned
	\renewenvironment{aligned}
		{\hbox to .88\displaywidth\bgroup$\displaystyle\originalaligned}
		{\endoriginalaligned$\hfil\egroup}
	\begin{gather}
		\begin{aligned}
			&\sum_{n\geq 1} a_n \left(\frac{3 (-1+3 n)}{5 n^2 (-1+2 n)^2}+\frac{(-2+7 n) \left(H_n-2 H_{2 n}\right)}{4 n^2 (-1+2 n)}\right) = \frac{7\zeta(3)}{40}
		\end{aligned} \\[3pt]
		\begin{aligned}
			&\sum_{n\geq 1} a_n \left(\frac{-11+28 n}{5 n^2 (-1+2 n)^2}+\frac{(-2+7 n) \left(H_n-H_{3 n}\right)}{n^2 (-1+2 n)}\right) = \frac{\zeta (3)}{10}
		\end{aligned} \\[3pt]
		\begin{aligned}
			&\sum_{n\geq 1} a_n \left(\frac{3 \left(14-73 n+88 n^2\right)}{5 n^2 (-1+2 n)^2 (-1+4 n)}+\frac{(-2+7 n) \left(H_n-4 H_{4 n}\right)}{n^2 (-1+2 n)}\right) = \frac{21 \zeta (3)}{5}
		\end{aligned} \\[3pt]
		\begin{aligned}
			&\sum_{n\geq 1} a_n \frac{3-18 n+28 n^2}{n^4 (-1+2 n)^3} = -\frac{\pi ^4}{45}
		\end{aligned} \\[3pt]
		\begin{aligned}
			&\sum_{n\geq 1} a_n\left(-\frac{9-58 n+102 n^2}{3 n^3 (-1+2 n)^3}+\frac{(-2+7 n) H_n^{(2)}}{n^2 (-1+2 n)} \right) = \frac{47 \pi ^4}{2160}
		\end{aligned} \\[3pt]
		\begin{aligned}
			&\sum_{n\geq 1} a_n \left(\frac{4 \left(6-43 n+84 n^2\right)}{3 n^3 (-1+2 n)^3}-\frac{\left(5-32 n+56 n^2\right) (H_n-H_{3 n})}{n^3 (-1+2 n)^2}\right) = -\frac{4 \pi ^4}{27}
		\end{aligned} \\[3pt]
		\begin{aligned}
			&\sum _{n\geq 1} a_n \left(\frac{2 \left(168 n^2-56 n+3\right)}{15 n^3 (2 n-1)^3}+\frac{\left(56 n^2-32 n+5\right) H_{3 n}}{5 n^3 (2 n-1)^2} -\frac{8 (7 n-2) H_{2 n}^{(2)}}{5 n^2 (2 n-1)}\right) =-\frac{37 \pi ^4}{1350}
		\end{aligned} \\[3pt]
		\begin{aligned}
			&\sum_{n\geq 1} a_n\left(-\frac{9-76 n+156 n^2}{3 n^3 (-1+2 n)^3}-\frac{\left(5-32 n+56 n^2\right) \left(H_n-2 H_{2 n}\right)}{2 n^3 (-1+2 n)^2}\right. \\
			&\left.\qquad +\frac{4 (-2+7 n) H_{2 n}^{(2)}}{n^2 (-1+2 n)}\right) = -\frac{11 \pi ^4}{270}
		\end{aligned} \\[3pt]
		\begin{aligned}
			&\sum_{n\geq 1} a_n\Bigg(-\frac{9-26 n+6 n^2}{24 n^3 (-1+2 n)^3}+\frac{3 (-1+3 n) (H_n-2 H_{2 n})}{n^2 (-1+2 n)^2} \\
			&\qquad +\frac{(-2+7 n) \left(5 H_n^2-20 H_n H_{2 n}+20 H_{2n}^2+4 H_{2 n}^{(2)}\right)}{8 n^2 (-1+2 n)} \Bigg) = -\frac{113 \pi ^4}{17280}
		\end{aligned} \\[3pt]
		\begin{aligned}
			&\sum_{n\geq 1} a_n \frac{17-136 n+408 n^2-640 n^3+560 n^4}{n^5 (-1+2 n)^4} = 180 \zeta (5)-\frac{56 \pi ^2 \zeta (3)}{3}
		\end{aligned} \\[3pt]
		\begin{aligned}
			&\sum_{n\geq 1} a_n \left(\frac{3 \left(17-102 n+123 n^2+107 n^3\right)}{17 n^4 (-1+2 n)^4}+\frac{(-2+7 n) H_n^{(3)}}{n^2 (-1+2 n)}\right) = \frac{360 \zeta (5)}{17}-\frac{155 \pi ^2 \zeta (3)}{68}
		\end{aligned} \\[3pt]
		\begin{aligned}
			&\sum_{n\geq 1} a_n \left(-\frac{12 (-1+3 n) \left(17-68 n+39 n^2\right)}{17 n^4 (-1+2 n)^4}+\frac{\left(5-32 n+56 n^2\right) H_n^{(2)}}{n^3 (-1+2 n)^2}\right) \\
			&\qquad = \frac{1103 \zeta (5)}{17}-\frac{322 \pi ^2 \zeta (3)}{51}
		\end{aligned} \\[3pt]
		\begin{aligned}
			&\sum _{n\geq 1} a_n \left(\frac{4 (3 n-1) (42 n-17)}{51 n^3 (2 n-1)^4}-\frac{\left(28 n^2-18 n+3\right) H_{2 n}}{3 n^4 (2 n-1)^3}\right) =\frac{7 \pi ^2 \zeta (3)}{51}-\frac{59 \zeta (5)}{51}
		\end{aligned} \\[3pt]
		\begin{aligned}
			&\sum _{n\geq 1} a_n \left(\frac{2 (662 n^3-266 n^2-34 n+17)}{51 n^4 (2 n-1)^4}-\frac{\left(28 n^2-18 n+3\right) H_n}{3 n^4 (2 n-1)^3}\right. \\
			&\left.\qquad +\frac{16 (7 n-2) H_{2 n}^{(3)}}{9 n^2 (2 n-1)}\right)=\frac{649 \zeta (5)}{51}-\frac{214 \pi ^2 \zeta (3)}{153}
		\end{aligned} \\[3pt]
		\begin{aligned}
			&\sum_{n\geq 1} a_n \left(-\frac{3 (-1+3 n) (17-68 n+22 n^2)}{17 n^4 (-1+2 n)^4}+\frac{(5-32 n+56 n^2) H_{2 n}^{(2)}}{n^3 (-1+2 n)^2}\right) \\
			&\qquad = \frac{815 \zeta (5)}{34}-\frac{140 \pi ^2 \zeta (3)}{51}
		\end{aligned} \\[3pt]
		\begin{aligned}
			&\sum_{n\geq 1}a_n\left(-\frac{8 (34-365 n+840 n^2)}{51 n^3 (-1+2 n)^4}+\frac{\left(3-18 n+28 n^2\right) \left(H_n-H_{3 n}\right)}{n^4 (-1+2 n)^3}\right) \\
			&\qquad = \frac{224 \pi ^2 \zeta (3)}{51}-\frac{1820 \zeta (5)}{51}
		\end{aligned} \\[3pt]
		\begin{aligned}
			&\sum_{n\geq 1}a_n\Bigg(-\frac{1683-7276 n-11908 n^2+62142 n^3}{765 n^4 (-1+2 n)^4}+\frac{(9-58 n+102 n^2) (H_n-H_{3 n})}{3 n^3 (-1+2 n)^3} \\
			&\qquad -\frac{(-11+28 n) H_n^{(2)}}{5 n^2 (-1+2 n)^2}-\frac{(-2+7 n) \left(H_n-H_{3 n}\right) H_n^{(2)}}{n^2 (-1+2 n)}\Bigg) = \frac{2408 \pi ^2 \zeta (3)}{765}-\frac{74639 \zeta (5)}{3060}
		\end{aligned} \\[3pt]
		\begin{aligned}
			&\sum_{n\geq 1} a_n\left(-\frac{-269+2978 n-13064 n^2+28432 n^3-30832 n^4+13504 n^5}{96 n^6 (-1+2 n)^5}\right. \\
			&\left.\qquad +\frac{(-2+7 n) H_n^{(4)}}{n^2 (-1+2 n)} \right) = \frac{7 \zeta (3)^2}{6}-\frac{19 \pi ^6}{15120}
		\end{aligned} \\[3pt]
		\begin{aligned}
			&\sum_{n\geq 1} a_n\left(-\frac{-31+406 n-2104 n^2+5360 n^3-6608 n^4+3008 n^5}{8 n^6 (-1+2 n)^5}\right. \\
			&\left.\qquad +\frac{\left(5-32 n+56 n^2\right) H_n^{(3)}}{n^3 (-1+2 n)^2} \right) = -2 \zeta (3)^2-\frac{\pi ^6}{756}
		\end{aligned} \\[3pt]
		\begin{aligned}
			&\sum _{n\geq 1} a_n \left(-\frac{3 (8 n^2-1)}{40 n^6 (2 n-1)^2}+\frac{(56 n^2-32 n+5) H_{2 n}^{(3)}}{5 n^3 (2 n-1)^2}\right) =-\frac{3 \zeta (3)^2}{5}-\frac{\pi ^6}{7560}
		\end{aligned} \\[3pt]
		\begin{aligned}
			&\sum_{n\geq 1} a_n\left( \frac{15-150 n+568 n^2-976 n^3+656 n^4}{4 n^6 (-1+2 n)^5} -\frac{\left(3-18 n+28 n^2\right) H_n^{(2)}}{n^4 (-1+2 n)^3}\right) = -\frac{\pi ^6}{378}
		\end{aligned} \\[3pt]
		\begin{aligned}
			&\sum _{n\geq 1} a_n \left(\frac{2240 n^5-8816 n^4+9968 n^3-5224 n^2+1330 n-133}{144 n^6 (2 n-1)^5}\right. \\
			&\left.\qquad +\frac{(28 n^2-18 n+3) H_{2 n}^{(2)}}{3 n^4 (2 n-1)^3}\right) =-\frac{5 \zeta (3)^2}{9}+\frac{\pi ^6}{1512}
		\end{aligned} \\[3pt]
		\begin{aligned}
			&\sum_{n\geq 1} a_n\Bigg(\frac{365-3362 n+10888 n^2-12432 n^3-3856 n^4+13760 n^5}{192 n^6 (-1+2 n)^5} \\
			&\qquad -\frac{\left(9-58 n+102 n^2\right) H_n^{(2)}}{3 n^3 (-1+2 n)^3}+\frac{(-2+7 n) (H_n^{(2)})^2}{2 n^2 (-1+2 n)} \Bigg) = \frac{7 \zeta (3)^2}{12}-\frac{\pi ^6}{324}
		\end{aligned} \\[3pt]
		\begin{aligned}
			&\sum_{n\geq 1} a_n\Bigg(-\frac{155-974 n+1592 n^2+1424 n^3-12272 n^4+22976 n^5}{96 n^6 (-1+2 n)^5} \\
			&\qquad -\frac{(17-136 n+408 n^2-640 n^3+560 n^4) \left(H_n-2 H_{2 n}\right)}{2 n^5 (-1+2 n)^4} \\
			&\qquad +\frac{32 (-2+7 n) H_{2 n}^{(4)}}{n^2 (-1+2 n)} \Bigg) = \frac{187 \pi ^6}{15120}-\frac{217 \zeta (3)^2}{6}
		\end{aligned} \\[3pt]
		\begin{aligned}
			&\sum_{n\geq 1} a_n \Bigg(\frac{-85+850 n-3400 n^2+6992 n^3-8240 n^4+5600 n^5}{3 n^6 (-1+2 n)^5} \\
			&\qquad -\frac{\left(17-136 n+408 n^2-640 n^3+560 n^4\right) \left(H_n-H_{3 n}\right)}{n^5 (-1+2 n)^4} \Bigg) = \frac{38 \pi ^6}{567}-\frac{392 \zeta (3)^2}{3}
		\end{aligned} \\[3pt]
		\begin{aligned}
			&\sum_{n\geq 1} a_n \left(\frac{-21+252 n-1196 n^2+2784 n^3-3120 n^4+1280 n^5}{8 n^7 (-1+2 n)^6}\right. \\
			&\left.\qquad -\frac{\left(3-18 n+28 n^2\right) H_n^{(3)}}{n^4 (-1+2 n)^3} \right) = \frac{\pi ^4 \zeta (3)}{45}
		\end{aligned} \\[3pt]
		\begin{aligned}
			&\sum_{n\geq 1} a_n\left( \frac{9 (-1+3 n)}{n^6 (-1+2 n)^2}-\frac{(5-32 n+56 n^2) (H_n^{(4)}-4 H_{2 n}^{(4)})}{n^3 (-1+2 n)^2}\right) \\
			&\qquad = -\frac{\pi ^4 \zeta (3)}{5}+\frac{62 \pi ^2 \zeta (5)}{3}-205 \zeta (7)
		\end{aligned} \\[3pt]
		\begin{aligned}
			&\sum_{n\geq 1} a_n \Bigg(\frac{20736 n^6-50560 n^5+51792 n^4-28032 n^3+8420 n^2-1332 n+87}{96 n^7 (2 n-1)^6} \\
			&\qquad -\frac{\left(28 n^2-18 n+3\right) H_{2 n}^{(3)}}{3 n^4 (2 n-1)^3}-\frac{\left(56 n^2-32 n+5\right) H_n^{(4)}}{4 n^3 (2 n-1)^2} \Bigg) \\
			&\qquad = -\frac{\pi ^4 \zeta (3)}{180}+\frac{31 \pi ^2 \zeta (5)}{18}-\frac{67 \zeta (7)}{4}
		\end{aligned} \\[3pt]
		\begin{aligned}
			&\sum_{n\geq 1} a_n \Bigg(-\frac{449-5964 n+33180 n^2-99488 n^3+170736 n^4-161152 n^5+67712 n^6}{8 n^7 (-1+2 n)^6} \\
			&\qquad +\frac{\left(17-136 n+408 n^2-640 n^3+560 n^4\right) H_n^{(2)}}{n^5 (-1+2 n)^4}+\frac{5 \left(5-32 n+56 n^2\right) H_n^{(4)}}{n^3 (-1+2 n)^2} \\
			&\qquad +\frac{4 (-2+7 n) H_n^{(5)}}{n^2 (-1+2 n)}\Bigg) = \frac{7 \pi ^4 \zeta (3)}{5}+\frac{185 \pi ^2 \zeta (5)}{3}-737 \zeta (7)
		\end{aligned} \\[3pt]
		\begin{aligned}
			&\sum_{n\geq 1} a_n \Bigg( \frac{-2+31 n-204 n^2+737 n^3-1573 n^4+1974 n^5-1343 n^6+385 n^7}{n^8 (-1+2 n)^7} \\
			&\qquad -\frac{\left(2-19 n+66 n^2-97 n^3+47 n^4\right) H_n^{(3)}}{n^5 (-1+2 n)^4} \\
			&\qquad +\frac{(-2+7 n)((H_n^{(3)})^2-H_n^{(6)})}{2 n^2 (-1+2 n)} \Bigg) = -\frac{ \pi ^2 \zeta (3)^2}{24}-\frac{31 \pi ^8}{725760}
		\end{aligned} \\[3pt]
		\begin{aligned}
			&\sum_{n\geq 1} a_n \Bigg(-\frac{63-882 n+5164 n^2-16232 n^3+29008 n^4-28128 n^5+11712 n^6}{8 n^8 (-1+2 n)^7} \\
			&\qquad +\frac{\left(15-150 n+568 n^2-976 n^3+656 n^4\right) H_n^{(2)}}{2 n^6 (-1+2 n)^5} \\
			&\qquad -\frac{(3-18 n+28 n^2) ((H_n^{(2)})^2-H_n^{(4)})}{n^4 (-1+2 n)^3} \Bigg) = \frac{\pi ^8}{1800}.
		\end{aligned}
	\end{gather}
	\endgroup
	For more formulas derivable in this manner, consult the accompanying Mathematica file. By extracting coefficients of order $\leq 6$, one can obtain a total of $1+4+12+29+59+108+180 = 393$ different formulas\footnote{This is obtained by summing the last row of the table.}. Most of them, especially those of high order, are too long to display here. \par
	Formulas (3.6), (3.10) and (3.11) were previously proved in \cite{hou2023taylor}; (3.15)-(3.21), (3.26), (3.32) imply seven conjectures of Sun \cite[(4.9)-(4.15)]{sun2024new}. A closely related conjecture \cite[(4.1)]{sun2024new}
	$$\sum_{n\geq 1} a_n\left(\frac{7n-2}{2n-1}H_{n-1}-\frac{1}{6n}\right) \stackrel{?}{=} \frac{\pi^2 \log 2}{6} - \frac{11\zeta(3)}{12}$$
	evidently does not follow from the formulas above because of the presence of $\log 2$. We encourage interested readers to seek other WZ-pairs that might prove it. \par
	
	A step-by-step computational guide for this example is available in the Mathematica package mentioned in the introduction. Readers who wish to see more formulas can consult the package.
\end{romexample}

\begin{romexample}\label{Ex_5}
	By applying Guillera's technique of a \textit{flawless WZ-pair} \cite[Example~4]{guillera2025wz} to the WZ-pair
	$$F(n,k) = \textsf{Seed9}(1 + a - 2 n, \frac12+ b - n, \frac12 + c - n, k + n),$$ or, equivalently, by taking $q\to 1$ in the $q$-analogue given in \cite[Example~XVII]{au2024wilfQ}, one can prove the following summation formula:
	{\small \begin{multline*}\sum_{n\geq 0}\frac{2^{-4 n} P(n) (\frac{1}{2}-b)_n (b+\frac{1}{2})_n (\frac{1}{2}-c)_n (c+\frac{1}{2})_n (-2 a+b+\frac{1}{2})_n (-2 a+c+\frac{1}{2})_n (-2 a+2 b+c+\frac{1}{2})_n (-2 a+b+2 c+\frac{1}{2})_n}{(1)_n (1-a)_n (-2 a+b+c+2 n+1) (-a+b+1)_n (-a+c+1)_n (-2 a+b+c+1)_n (-2 a+b+c+1)_{2 n} (-a+b+c+1)_n}  \\
		= \frac{\Gamma (1-a) 2^{-4 a+2 b+2 c+5} \Gamma (-a+b+1) \Gamma (-a+c+1) \Gamma (-2 a+b+c+1)^2 \Gamma (-a+b+c+1)}{\pi  \Gamma \left(-2 a+b+\frac{1}{2}\right) \Gamma \left(-2 a+c+\frac{1}{2}\right) \Gamma \left(-2 a+2 b+c+\frac{1}{2}\right) \Gamma \left(-2 a+b+2 c+\frac{1}{2}\right)} .\end{multline*} }
	where $P(n)$ is the following polynomial:
	{\small $
			P(n) = 8 n^2 (240 a^2-120 a (2 b+2 c+1)+54 b^2+60 b (2 c+1)+54 c^2+60 c+13)+(4 a-2 b-1) (4 a-2 c-1) (16 a^2-8 a (3 b+3 c+1)+8 b^2+b (20 c+6)+8 c^2+6 c+1)+8 n (-144 a^3+8 a^2 (27 b+27 c+14)-2 a \left(48 b^2+4 b (27 c+14)+48 c^2+56 c+13\right)+12 b^3+24 b^2 (2 c+1)+b \left(48 c^2+56 c+13\right)+12 c^3+24 c^2+13 c+2)+32 n^3 (-42 a+21 b+21 c+10)+336 n^4. $ }
Restricting the formula above to $a=b=c=0$ recovers the famous formula (conjectured by Gourevich and proved by the author \cite{au2025wilf}): 	$$\sum_{n\geq 0} \left(\frac{1}{2^6}\right)^n \frac{(\frac12)_n^7}{(1)_n^7} (168 n^3+76 n^2+14 n+1) = \frac{32}{\pi^3}.$$
Let $S := S(a,b,c;n)$ denote the summand on the left-hand side of the formula above. 

\begin{conjecture}
	The generating function
	$$\sum_{N\geq 0} \dim V_N(S) t^N \stackrel{?}{=} \frac{1}{(1-t)(1-t^2)(1-t^4)} = 1+t+2t^2+2t^3+4t^4+4t^5+6t^6+\cdots.$$
\end{conjecture}

We list some examples of identities in $V_N(S)$, including only the simplest-looking ones. To abbreviate the notation, write
$$\begin{aligned}
&a_n = \left(\frac{1}{2^6}\right)^n \frac{(\frac12)_n^7}{(1)_n^7} \\
&p_0 = 1+8 n+28 n^2,\quad p_1 = 7+76 n+252 n^2,\quad p_2 = 1+6n, \quad p_3 = 83 + 546 n.
\end{aligned}$$
Then we obtain
	\begin{align*}
		&\color{NavyBlue}\sum_{n\geq 0} a_n\left(\frac{p_1}{7}+\left(-H_n+H_{2 n}\right) p_0 p_2\right) = \frac{320 \log (2)}{7 \pi ^3} \\
		&\color{NavyBlue}\sum_{n\geq 0} a_n p_2 \left(\frac{1}{336}+\frac{\left(-5 H_n^{(2)}+16 H_{2 n}^{(2)}\right) p_0}{2688}\right) = \frac{5}{504 \pi } \\
		&\sum_{n\geq 0} a_n \left(\left(-H_n+H_{2 n}\right) p_1+\frac{1}{32} \left(112 \left(H_n\right){}^2-224 H_n H_{2 n}+112 \left(H_{2 n}\right){}^2+3 H_n^{(2)}\right) p_0 p_2+\frac{p_3}{28}\right) = \frac{1600 \log ^2(2)}{7 \pi ^3}-\frac{23}{14 \pi } \\
		&\color{NavyBlue}\sum_{n\geq 0} a_n \left(\frac{1}{1+2 n}+\frac{1}{64} \left(-7 H_n^{(4)}+128 H_{2 n}^{(4)}\right) p_0 p_2 \right) = \frac{61 \pi }{180} \\
		&\sum_{n\geq 0} a_n\left(\frac{1}{(1+2 n)}+\frac{-5 H_n^{(2)}+16 H_{2 n}^{(2)}}{4} p_2+\frac{\left(50 \left(H_n^{(2)}\right){}^2+3 H_n^{(4)}-320 H_n^{(2)} H_{2 n}^{(2)}+512 \left(H_{2 n}^{(2)}\right){}^2\right) p_0}{128} p_2 \right) = \frac{16\pi }{45} \\
		&\sum_{n\geq 0} a_n\Bigg(-\frac{1+7 (1+2 n) (H_n-H_{2n})}{(1+2 n)^2}+\frac{1}{64} \left(-7 H_n^{(4)}+128 H_{2 n}^{(4)}\right) p_1 \\ &\quad +\frac{1}{64} \left(49 H_n H_n^{(4)}-49 H_{2 n} H_n^{(4)}+14 H_n^{(5)}-896 H_n H_{2 n}^{(4)}+896 H_{2 n} H_{2 n}^{(4)}-512 H_{2 n}^{(5)}\right) p_0 p_2 \Bigg) = \frac{61 \pi  \log (2) }{18} -\frac{249 \zeta (5)}{\pi ^3}
	\end{align*}
	
	The identities involving first-, second-, and fourth-order harmonic numbers were conjectured by Sun \cite[(181),(183),(185)]{sun2022conjectures}; the identity involving third-order harmonic numbers is more complicated (the fourth equation). \par
	A step-by-step computational guide for this example is available in the Mathematica package mentioned in the introduction. 
\end{romexample}
\section{More Examples}
For the WZ-pairs in this section, at least one of the following terms in Proposition \ref{WZ_prop} will be nonzero: $$\sum_{n\geq 0} g(n)\qquad \text{and} \qquad \lim_{n\to \infty} \sum_{k\geq 0} F(n,k).$$
There are several ways to address these difficulties, including the following:
\begin{itemize}[leftmargin=*]
	\item Au's technique \cite{au2025wilf, au2025multiple}, based on asymptotic analysis, which determines $\lim_{n\to \infty} \sum_{k\geq 0} F(n,k)$ explicitly. 
	\item Guillera's technique \cite{guillera2025wz} of \textit{flawless WZ-pairs}; it works well for $1/\pi^k$ formulas but less well for series leading to odd zeta values (e.g., Example \ref{Ex_10}). 
\end{itemize}
We will not, however, address these issues in this article. Consequently, some of the identities listed below are, strictly speaking, conjectural because explicit expressions for $g(n)$ and $\lim_{n\to \infty} \sum_{k\geq 0} F(n,k)$ are unavailable. They are indicated by an question mark above the equality sign $\stackrel{?}{=}$.

\begin{romexample}\label{Ex_6}
	This example concerns a harmonic extension of the series
	\begin{equation}\label{aux_4}\sum_{n\geq 1} \left(\frac{1}{4}\right)^n \frac{(1)_n^3}{(\frac12)_n^3} \frac{3n-1}{n^3} = \frac{\pi^2}{2}.\end{equation}
	There are two WZ-seeds that give this sum, namely,
   $$\begin{aligned}F_1(n,k) &= \textsf{Seed2}\left(a-n+1,b+\frac{1}{2},c-n+2,k+d+n-\frac{1}{2}\right), \\ F_2(n,k) &= \textsf{Gauss2F1}\left(a-n+1,b-n+1,c+1,k+d+n-\frac{1}{2}\right).\end{aligned}$$
   From these pairs, we can write down the relevant summands below\footnote{More precisely, they are obtained from $G_i(n,0)$ after omitting factors independent of $n$.}.
	$$S_1(a,b,c,d;n) := \frac{2^{-2 n} (a+2 b-2 c+2 d+3 n-1) (1-a)_{n-1} (a+2 d+1)_{n-1} (b+d+1)_{n-1}  (a+2 b-2 c+1)_{n-1}}{\left(d+\frac{1}{2}\right)_n (a-2 c+1)_{n-1} \left(b-c+\frac{1}{2}\right)_n \left(a+b-c+d+\frac{1}{2}\right)_n}$$
	
	$$S_2(a,b,c,d;n) := \frac{P(n) (1-a)_{n-1} (1-b)_{n-1} (-a+c+1)_{n-1} (-b+c+1)_{n-1} }{\left(d+\frac{1}{2}\right)_n \left(c+d+\frac{1}{2}\right)_n (-a-b+c+1)_{2 n}},$$
	where $P(n) =(b - c) (-1 + 2 c + 2 d) + 
	a (-1 - 2 b + 2 c + 2 d) + (2 + 4 a + 4 b - 6 c - 4 d) n - 6 n^2$. The table shows that each summand provides new relations.
	
	\begin{table}[h]
	\begin{tabular}{|c|c|c|c|c|c|c|c|}
		\hline
		$N$                          & $0$ & $1$ & $2$  & $3$  & $4$  & $5$  & $6$   \\ \hline
		$\dim V_N(S_1)$              & $1$ & $2$ & $4$ & $7$ & $11$ & $16$ & $23$ \\ \hline
		$\dim V_N(S_2)$              & $1$ & $2$ & $5$  & $8$ & $14$ & $20$ & $30$  \\ \hline
		$\dim (V_N(S_1) + V_N(S_2))$ & $1$ & $2$ & $5$ & $9$ & $16$ & $24$ & $37$ \\ \hline
	\end{tabular}
		\caption{\small Number of independent summation formulas that can be generated by $S_1(a,b,c,d;n)$ and $S_2(a,b,c,d;n)$, separately and in combination.}
	\end{table}
	\begin{conjecture}The generating functions of $\dim V_N(S_1)$, $\dim V_N(S_2)$, and $\dim (V_N(S_1) + V_N(S_2))$ are rational functions. More precisely, 
	$$\sum_{N\geq 0} \dim V_N(S_1) t^N \stackrel{?}{=} \frac{1}{(1-t)^2(1-t^2)(1-t^3)},\qquad \sum_{N\geq 0} \dim V_N(S_2) t^N\stackrel{?}{=} \frac{1}{(1-t)^2 (1-t^2)^2}.$$
	\end{conjecture}
	Write $a_n = \left(\frac{1}{4}\right)^n \frac{(1)_n^3}{(\frac12)_n^3}$. We list some examples of identities in $V_N(S_1)+V_N(S_2)$, including only the simplest-looking ones. 
	\begingroup
	\let\originalaligned\aligned
	\let\endoriginalaligned\endaligned
	\renewenvironment{aligned}
		{\hbox to .88\displaywidth\bgroup$\displaystyle\originalaligned}
		{\endoriginalaligned$\hfil\egroup}
	\begin{gather}
		\begin{aligned}
			&\sum_{n\geq 1} a_n \frac{(-1+3 n) \left(-1+n H_n\right)}{n^4} = -\frac{2\pi^2}{3} \log (2)+\frac{14 \zeta (3)}{3}
		\end{aligned} \\[3pt]
		\begin{aligned}
			&\sum_{n\geq 1} a_n \left(\frac{-1+4 n}{2 n^4}-\frac{(-1+3 n) H_{2 n}}{n^3}\right) = \frac{1}{3} \pi ^2 \log (2)-\frac{35 \zeta (3)}{6}
		\end{aligned} \\[3pt]
		\begin{aligned}
			&\sum_{n\geq 1} a_n \left(\frac{-2+5 n}{2 n^5}-\frac{(-1+3 n) H_n^{(2)}}{n^3}\right) = -\frac{\pi ^4}{48}
		\end{aligned} \\[3pt]
		\begin{aligned}
			&\sum_{n\geq 1} a_n \left(\frac{-2+n}{8 n^5}-\frac{(-1+3 n) H_{2 n}^{(2)}}{n^3} \right) = -\frac{13 \pi ^4}{192}
		\end{aligned} \\[3pt]
		\begin{aligned}
			&\sum_{n\geq 1} a_n \left(-\frac{-6+17 n}{6 n^5}-\frac{(-1+3 n) H_n \left(-2+n H_n\right)}{n^4} \right) \\
			&\qquad \stackrel{?}{=} -\frac{128 \text{Li}_4\left(\frac{1}{2}\right)}{9}+\frac{533 \pi ^4}{6480} -\frac{16}{27} \log ^4(2)-\frac{8}{27} \pi ^2 \log ^2(2)
		\end{aligned} \\[3pt]
		\begin{aligned}
			&\sum_{n\geq 1} a_n \left(-\frac{-2+7 n}{4 n^5}+\frac{(-1+4 n) H_n}{2 n^4}-\frac{(-1+3 n) \left(-1+n H_n\right) H_{2 n}}{n^4} \right) \\
			&\qquad \stackrel{?}{=} -\frac{112 \text{Li}_4\left(\frac{1}{2}\right)}{9}+\frac{629 \pi ^4}{12960}-\frac{14}{27} \log ^4(2)+\frac{2}{27} \pi ^2 \log ^2(2)
		\end{aligned} \\[3pt]
		\begin{aligned}
			&\sum_{n\geq 1} a_n \left(-\frac{-6+19 n}{24 n^5}+\frac{(-1+4 n) H_{2 n}}{n^4}-\frac{(-1+3 n) (H_{2 n})^2}{n^3}\right) \\
			&\qquad \stackrel{?}{=} -\frac{80 \text{Li}_4\left(\frac{1}{2}\right)}{9}-\frac{17 \pi ^4}{5184}-\frac{10}{27} \log ^4(2)+\frac{4}{27} \pi ^2 \log ^2(2)
		\end{aligned} \\[3pt]
		\begin{aligned}
			&\sum_{n\geq 1} a_n \left(-\frac{(-1+3 n) \left(-8+7 n^3 H_n^{(3)}+8 n^3 H_{2 n}^{(3)}\right)}{8 n^6}\right) = -\frac{\pi ^2 \zeta (3)}{2}
		\end{aligned} \\[3pt]
		\begin{aligned}
			&\sum_{n\geq 1} a_n \left(-\frac{-5+13 n}{3 n^6}+\frac{(-2+5 n) H_n}{2 n^5}-\frac{(-1+3 n) \left(-3 H_n^{(2)}+3 n H_n H_n^{(2)}-2 n H_n^{(3)}\right)}{3 n^4}\right) \\
			&\qquad \stackrel{?}{=} -\frac{17 \pi ^2 \zeta (3)}{12}+\frac{341 \zeta (5)}{24}+\frac{1}{36} \pi ^4 \log (2)
		\end{aligned} \\[3pt]
		\begin{aligned}
			&\sum_{n\geq 1} a_n \left(-\frac{-13+29 n}{6 n^6}+\frac{(-2+n) H_n}{4 n^5}+\frac{(-1+3 n) \left(5 n H_n^{(3)}+6 H_{2 n}^{(2)}-6 n H_n H_{2 n}^{(2)}\right)}{3 n^4} \right) \\
			&\qquad \stackrel{?}{=} -\frac{19 \pi ^2 \zeta (3)}{8}+\frac{713 \zeta (5)}{48}+\frac{13}{72} \pi ^4 \log (2)
		\end{aligned} \\[3pt]
		\begin{aligned}
			&\sum_{n\geq 1} a_n \left(\frac{-7+19 n}{3 n^6}-\frac{(-6+17 n) H_n}{2 n^5}-\frac{(-1+3 n) \left(-9 (H_n)^2+3 n (H_n)^3+4 n H_n^{(3)}\right)}{3 n^4} \right) \\
			&\qquad \stackrel{?}{=} -\frac{512 \text{Li}_5\left(\frac{1}{2}\right)}{9}+\frac{67 \pi ^2 \zeta (3)}{36}+\frac{1705 \zeta (5)}{72} +\frac{64 \log ^5(2)}{135}+\frac{32}{81} \pi ^2 \log ^3(2)-\frac{533 \pi ^4 \log (2)}{1620}
		\end{aligned} \\[3pt]
		\begin{aligned}
			&\sum_{n\geq 1} a_n\left(\frac{-71+200 n}{12 n^6}-\frac{(-6+19 n) H_{2 n}}{4 n^5} +\frac{3 (-1+4 n) \left(H_{2 n}\right){}^2}{n^4}-\frac{(-1+3 n) (6H_{2 n}^3+17 H_n^{(3)})}{3 n^3}\right) \\
			&\qquad \stackrel{?}{=} -\frac{320 \text{Li}_5\left(\frac{1}{2}\right)}{9}+\frac{545 \pi ^2 \zeta (3)}{72}-\frac{24149 \zeta (5)}{288} +\frac{8 \log ^5(2)}{27}-\frac{16}{81} \pi ^2 \log ^3(2)+\frac{17 \pi ^4 \log (2)}{1296}
		\end{aligned} \\[3pt]
		\begin{aligned}
			&\sum_{n\geq 1} a_n\left(-\frac{-3+7 n}{2 n^7}+\frac{(-2+5 n) H_n^{(2)}}{n^5}-\frac{(-1+3 n) \left(2 \left(H_n^{(2)}\right){}^2-H_n^{(4)}\right)}{2 n^3}\right) \stackrel{?}{=} -\frac{\pi ^6}{1440}
		\end{aligned} \\[3pt]
		\begin{aligned}
			&\sum_{n\geq 1} a_n\left(\frac{-1+9 n}{4 n^7}+\frac{(-2+n) H_n^{(2)}}{4 n^5} \right. \\
			&\left.\qquad +\frac{(-2+5 n) H_{2 n}^{(2)}}{n^5}-\frac{(-1+3 n) \left(3 H_n^{(4)}+8 H_n^{(2)} H_{2 n}^{(2)}\right)}{4 n^3} \right) \stackrel{?}{=} -\frac{17 \pi ^6}{2880}
		\end{aligned} \\[3pt]
		\begin{aligned}
			&\sum_{n\geq 1} a_n\left(\frac{-21+73 n}{8 n^7}+\frac{(-2+n) H_{2 n}^{(2)}}{n^5} \right. \\
			&\left.\qquad -\frac{(-1+3 n) \left(21 H_n^{(4)}+32 \left(H_{2 n}^{(2)}\right){}^2+32 H_{2 n}^{(4)}\right)}{8 n^3}\right) \stackrel{?}{=} -\frac{337 \pi ^6}{5760}
		\end{aligned} \\[3pt]
		\begin{aligned}
			&\sum_{n\geq 1} a_n\left(-\frac{(-1+3 n) \left(-1+n H_n\right) \left(-8+7 n^3 H_n^{(3)}+8 n^3 H_{2 n}^{(3)}\right)}{8 n^7}\right) \\
			&\qquad \stackrel{?}{=} -\frac{7 \zeta (3)^2}{3}+\frac{2}{3} \pi ^2 \zeta (3) \log (2)-\frac{\pi ^6}{224}
		\end{aligned} \\[3pt]
		\begin{aligned}
			&\sum_{n\geq 1} a_n\Bigg(\frac{-129+299 n}{36 n^7}-\frac{(-13+29 n) H_n}{3 n^6}+\frac{(-2+n) (H_n)^2}{4 n^5}-\frac{(-6+17 n) H_{2 n}^{(2)}}{3 n^5} \\
			&\qquad +\frac{(-1+3 n) \left(-40 H_n^{(3)}+40 n H_n H_n^{(3)}+3 n H_n^{(4)}+48 H_n H_{2 n}^{(2)}-24 n (H_n)^2 H_{2 n}^{(2)}\right)}{12 n^4} \Bigg) \\
			&\qquad \stackrel{?}{=} \frac{184}{9}\zeta(\overline{5},1)-\frac{179 \zeta (3)^2}{6}+\frac{19}{3} \pi ^2 \zeta (3) \log (2) -\frac{713}{18} \zeta (5) \log (2)+\frac{827 \pi ^6}{25920}-\frac{13}{54} \pi ^4 \log ^2(2)
		\end{aligned}
	\end{gather}
	\endgroup
	Among the formulas above, those not marked with $\stackrel{?}{=}$ were conjectured by Sun \cite{sun2022conjectures} and proved in \cite{li2023infinite4}.
\end{romexample}

\begin{romexample}\label{Ex_7}
This example considers a harmonic extension of Ramanujan's $1/\pi$ formula \cite{cohen2021rational, chudnovsky1988approximations}:
$$\sum_{n\geq 0} \left(\frac{1}{2^6}\right)^n \frac{(\frac12)_n^3}{(1)_n^3} (42n+5) = \frac{16}{\pi}.$$
$$\begin{aligned}
F_1(n,k) = \textsf{Seed2}(a-n+\frac{1}{2},b+\frac{1}{2},c+n+1,d+k+n),\\
F_2(n,k) = \textsf{Dixon6F5}(a-n+\frac{1}{2},b+1,c+1,d+k+n+\frac{1}{2}),\\
F_3(n,k) = \textsf{Seed9}(a-2 n+1,b-n+\frac{1}{2},c-n+\frac{1}{2},k+d+n-\frac{1}{6}).\\
\end{aligned}$$

As in the previous examples, we define $S_i$ from the WZ-mate of $F_i(n,k)$ by retaining only the factors that depend on $n$.\footnote{For the second WZ-seed, one has to do an additional telescoping for $G_2(n,k)$ to obtain the $S_2$ provided in the Mathematica notebook}	

\begin{table}[h]
	\begin{tabular}{|c|c|c|c|c|c|c|c|}
		\hline
		$N$                          & $0$ & $1$ & $2$  & $3$  & $4$  & $5$  & $6$    \\ \hline
		$\dim V_N(S_1)$              & $1$ & $2$ & $6$ & $11$ & $20$ & $31$ & $47$  \\ \hline
		$\dim V_N(S_2)$              & $1$ & $1$ & $3$ & $3$ & $7$ & $7$ & $13$  \\ \hline
		$\dim V_N(S_3)$              & $1$ & $2$ & $5$ & $8$ & $14$ & $20$ & $30$  \\ \hline
		$\dim (V_N(S_1) + V_N(S_2))$              & $1$ & $2$ & $6$ & $11$ & $21$ & $32$ & $51$  \\ \hline
		$\dim (V_N(S_1) + V_N(S_2)+V_N(S_3))$              & $1$ & $3$ & $8$ & $16$ & $30$ & $47$ & $74$  \\ \hline
	\end{tabular}
	\caption{\small Number of independent summation formulas that can be generated by $S_i(a,b,c,d;n)$.}
\end{table}

\begin{conjecture}
	The generating function for each row of the table above is a rational function. More precisely,
	\begin{align*}
		\sum_{N= 0}^\infty \dim V_N(S_1) t^N &\stackrel{?}{=}\frac{1-t+2t^2-t^3}{(1-t)^3(1-t^2)} \\
		\sum_{N= 0}^\infty \dim V_N(S_2) t^N &\stackrel{?}{=} \frac{1}{(1-t)(1-t^2)^2(1-t^4)} \\
		\sum_{N= 0}^\infty \dim V_N(S_3) t^N &\stackrel{?}{=} \frac{1}{(1-t)^2(1-t^2)^2}\end{align*}
\end{conjecture}

When extracting power-series coefficients from $S_1$ and $S_2$, the only harmonic numbers that appear are $H_{dn}^{(i)}, d\in \{1,2,3,6\}$; for $S_3$, we also require the new harmonic number $U_n^{(k)} := \sum_{i=0}^{n-1} (i+\frac16)^{-k}$. We list some examples of identities in $V_N(S_1)+V_N(S_2) + V_N(S_3)$, including only the simplest-looking ones. To abbreviate the notation, set $$a_n :=  \left(\frac{1}{2^6}\right)^n \frac{(\frac12)_n^3}{(1)_n^3}, \qquad p_0 := 42n+5, \qquad p_2 = 35+162n,\qquad p_5 = 77+366n.$$ Then
\begingroup
\let\originalaligned\aligned
\let\endoriginalaligned\endaligned
\renewenvironment{aligned}
	{\hbox to .88\displaywidth\bgroup$\displaystyle\originalaligned}
	{\endoriginalaligned$\hfil\egroup}
\begin{gather}
	\begin{aligned}
		&\sum_{n\geq 0} a_n\left(1+\frac{1}{7} \left(-H_n+H_{2 n}\right) (5+42 n)\right) = \frac{32 \log (2)}{7 \pi }
	\end{aligned} \\[3pt]
	\begin{aligned}
		&\sum_{n\geq 0} a_n\left(\frac{31+162 n}{3+18 n}+\left(-\frac{H_n}{2}-H_{3 n}+2 H_{6 n}\right) (5+42 n)\right) \stackrel{?}{=} \frac{48 \log (2)}{\pi }
	\end{aligned} \\[3pt]
	\begin{aligned}
		&\sum_{n\geq 0} a_n \left(\frac{35+162 n}{1+6 n}+\left(-\frac{H_n}{2}+2 U_n^{(1)}\right) (5+42 n) \right) \stackrel{?}{=} \frac{32}{\sqrt{3}}+\frac{80 \log (2)}{\pi }
	\end{aligned} \\[3pt]
	\begin{aligned}
		&\sum_{n\geq 0} a_n \left( \frac{1}{2} (42 n+5) H_n^{(2)}+\frac{46}{2 n+1} \right) = \frac{44\pi}{3}
	\end{aligned} \\[3pt]
	\begin{aligned}
		&\sum_{n\geq 0} a_n \left( (42 n+5) H_{2 n}^{(2)}+\frac{25}{2 n+1} \right) = 8\pi
	\end{aligned} \\[3pt]
	\begin{aligned}
		&\sum_{n\geq 0} a_n \left(\frac{43+396 n+972 n^2}{9 (1+2 n) (1+6 n)^2}+\left(-\frac{H_{3 n}^{(2)}}{4}+H_{6 n}^{(2)}\right) (5+42 n) \right) \stackrel{?}{=} \frac{14 \pi }{9}
	\end{aligned} \\[3pt]
	\begin{aligned}
		&\sum_{n\geq 0} a_n \left( -2 H_n+2 H_{2 n}-\frac{1}{2 n+1} + \frac{1}{7} (42 n+5) (H_n-H_{2 n})^2 \right) = \frac{64 \log ^2(2)}{7 \pi }-\frac{16 \pi }{21}
	\end{aligned} \\[3pt]
	\begin{aligned}
		&\sum_{n\geq 0} a_n \left(\frac{(6 n+5) (78 n+17)}{9 (2 n+1) (6 n+1)^2} + \frac{p_0}{4} (H_n+2 H_{3 n}-4 H_{6 n})^2\right. \\
		&\left.\qquad -\frac{(162 n+31) \left(H_n+2 H_{3 n}-4 H_{6 n}\right)}{18 n+3} \right) \stackrel{?}{=} \frac{144 \log ^2(2)}{\pi }-\frac{34 \pi }{9}
	\end{aligned} \\[3pt]
	\begin{aligned}
		&\sum_{n\geq 0} a_n\Bigg(14 U_n^{(1)}-\frac{8}{(6 n+1)^2}-\frac{23}{6 n+3}+\frac{1}{6} \left(3 H_n^2+12 H_{2 n} U_n^{(1)}-3 H_n \left(H_{2 n}+4 U_n^{(1)}\right)-2 U_n^{(2)}\right) p_0 \\
		&\qquad +\frac{H_{2 n} p_2}{1+6 n}-\frac{H_n p_5}{2+12 n}\Bigg) \stackrel{?}{=} -\frac{120 L_{-3}(2)}{\pi }-\frac{104 \pi }{9}+\frac{160 \log ^2(2)}{\pi }+\frac{64 \log (2)}{\sqrt{3}}
	\end{aligned} \\[3pt]
	\begin{aligned}
		&\sum_{n\geq 0} a_n \left(-\frac{9}{(1+2 n)^2}+\frac{1}{3} \left(-2 H_n^{(3)}+17 H_{2 n}^{(3)}\right) p_0\right) \stackrel{?}{=} \frac{80 \zeta (3)}{\pi }-\frac{128 G}{3}
	\end{aligned} \\[3pt]
	\begin{aligned}
		&\sum_{n\geq 0} a_n \Bigg( \frac{-479-1275 (1+2 n) H_n+1275 (1+2 n) H_{2 n}+357 H_{2 n}^{(2)}+1428 n H_{2 n}^{(2)}+1428 n^2 H_{2 n}^{(2)}}{(1+2 n)^2} \\
		&\qquad +\left(-4 H_n^{(3)}+51 \left(-H_n+H_{2 n}\right) H_{2 n}^{(2)}\right) p_0\Bigg) \stackrel{?}{=} -2432 G-\frac{64 \zeta (3)}{\pi }+816 \pi \log (2)
	\end{aligned} \\[3pt]
	\begin{aligned}
		&\sum_{n\geq 0} a_n \left(\frac{25}{9 (1+2 n)^3}+\left(-\frac{H_n^{(4)}}{16}+H_{2 n}^{(4)}\right) (5+42 n) \right) \stackrel{?}{=} \frac{5 \pi ^3}{54}
	\end{aligned} \\[3pt]
	\begin{aligned}
		&\sum_{n\geq 0} a_n \left(\frac{4 \left(-5620+621 (1+2 n)^2 H_n^{(2)}\right)}{27 (1+2 n)^3}+\frac{1}{2} \left(\left(H_n^{(2)}\right){}^2+H_n^{(4)}\right) (5+42 n)\right) \stackrel{?}{=} -\frac{10874 \pi ^3}{405}
	\end{aligned} \\[3pt]
	\begin{aligned}
		&\sum_{n\geq 0} a_n \left(\frac{25 H_n^{(2)}}{1+2 n}+\frac{4 \left(-3103+621 (1+2 n)^2 H_{2 n}^{(2)}\right)}{27 (1+2 n)^3}\right. \\
		&\left.\qquad +\left(\frac{H_n^{(4)}}{4}+H_n^{(2)} H_{2 n}^{(2)}\right) (5+42 n) \right) \stackrel{?}{=} -\frac{6002 \pi ^3}{405}
	\end{aligned} \\[3pt]
	\begin{aligned}
		&\sum_{n\geq 0} a_n \left( \frac{5 \left(-673+270 (1+2 n)^2 H_{2 n}^{(2)}\right)}{27 (1+2 n)^3}+\left(\frac{H_n^{(4)}}{16}+\left(H_{2 n}^{(2)}\right){}^2\right) (5+42 n) \right) \stackrel{?}{=} -\frac{3251 \pi ^3}{810}
	\end{aligned} \\[3pt]
	\begin{aligned}
		&\sum_{n\geq 0} a_n \Bigg(9 p_0 \left(16 \left(H_n-H_{2 n}\right) \left(2 H_n^{(3)}-17 H_{2 n}^{(3)}\right)-H_n^{(4)}\right) \\
		&\qquad +\frac{3888 \left(H_n-H_{2 n}\right)}{(2 n+1)^2}-2016 H_n^{(3)}+17136 H_{2 n}^{(3)}+\frac{9392}{(2 n+1)^3}\Bigg) \\
		&\qquad \stackrel{?}{=} 73728 \Im\left(\text{Li}_3\left(\frac{1}{2}+\frac{i}{2}\right)\right)+\frac{69120 \zeta (3) \log (2)}{\pi }-\frac{22904 \pi ^3}{15}-2304 \pi  \log ^2(2)
	\end{aligned} \\[3pt]
	\begin{aligned}
		&\sum_{n\geq 0} a_n \Bigg(H_n^{(2)} \left((17 H_{2 n}^{(3)}-2 H_n^{(3)})p_0-\frac{27}{(2 n+1)^2}\right) \\
		&\qquad +4 \left(5 p_0 H_n^{(5)}+\frac{391 H_{2 n}^{(3)}-46 H_n^{(3)}}{2 n+1}-159 p_0 H_{2 n}^{(5)}-\frac{4}{(2 n+1)^4}\right) \Bigg) \\
		&\qquad \stackrel{?}{=} 1536 L_{-4}(4)+440 \pi \zeta (3)-\frac{9856 \zeta (5)}{\pi }
	\end{aligned} \\[3pt]
	\begin{aligned}
		&\sum_{n\geq 0} a_n \Bigg(\left(H_{2 n}-H_n\right) \left((144 H_{2 n}^{(4)}-9H_n^{(4)})p_0+\frac{400}{(2 n+1)^3}\right) \\
		&\qquad -63 H_n^{(4)}+6 (42 n+5) H_n^{(5)}-8064 n H_{2 n}^{(5)} \\
		&\qquad +16 \left(63 H_{2 n}^{(4)}-60 H_{2 n}^{(5)}-\frac{25}{(2 n+1)^4}\right) \Bigg) \stackrel{?}{=} \frac{80}{3} \pi ^3 \log (2)-\frac{2976 \zeta (5)}{\pi }
	\end{aligned} \\[3pt]
	\begin{aligned}
		&\sum_{n\geq 0} a_n \Bigg(-\frac{50 H_n^{(3)}}{2 n+1}+\frac{425 H_{2 n}^{(3)}}{2 n+1}+p_0 \left(6 H_n^{(5)}-2 H_n^{(3)} H_{2 n}^{(2)}+17 H_{2 n}^{(2)} H_{2 n}^{(3)}-191 H_{2 n}^{(5)}\right) \\
		&\qquad -\frac{27 H_{2 n}^{(2)}}{(2 n+1)^2}-\frac{2}{(2 n+1)^4} \Bigg) \stackrel{?}{=} 512 L_{-4}(4)+120 \pi \zeta (3)-\frac{2960 \zeta (5)}{\pi }
	\end{aligned} \\[3pt]
	\begin{aligned}
		&\sum_{n\geq 0} a_n \Bigg(p_0 \left(-27 H_n^{(2)} H_n^{(4)}-172 H_n^{(6)}+432 H_n^{(2)} H_{2 n}^{(4)}+11008 H_{2 n}^{(6)}\right) \\
		&\qquad +\frac{1200 H_n^{(2)}}{(2 n+1)^3}-\frac{2484 \left(H_n^{(4)}-16 H_{2 n}^{(4)}\right)}{2 n+1}-\frac{192}{(2 n+1)^5} \Bigg) \stackrel{?}{=} \frac{44 \pi ^5}{15}
	\end{aligned} \\[3pt]
	\begin{aligned}
		&\sum_{n\geq 0} a_n \Bigg(p_0 \left(27 \left(H_n^{(2)}\right){}^3+81 H_n^{(4)} H_n^{(2)}+13022 H_n^{(6)}-829952 H_{2 n}^{(6)}\right) \\
		&\qquad +\frac{7452 \left(H_n^{(2)}\right){}^2}{2 n+1}-\frac{134880 H_n^{(2)}}{(2 n+1)^3}+\frac{7452 H_n^{(4)}}{2 n+1}+\frac{428928}{(2 n+1)^5} \Bigg) \stackrel{?}{=} \frac{120706 \pi ^5}{105}
	\end{aligned} \\[3pt]
	\begin{aligned}
		&\sum_{n\geq 0} a_n \Bigg(p_0 \left(-155 H_n^{(6)}-108 H_n^{(4)} H_{2 n}^{(2)}+1728 H_{2 n}^{(2)} H_{2 n}^{(4)}+9920 H_{2 n}^{(6)}\right) \\
		&\qquad -\frac{2700 H_n^{(4)}}{2 n+1}+\frac{4800 H_{2 n}^{(2)}}{(2 n+1)^3}+\frac{43200 H_{2 n}^{(4)}}{2 n+1}-\frac{1152}{(2 n+1)^5} \Bigg) \stackrel{?}{=} 0
	\end{aligned} \\[3pt]
	\begin{aligned}
		&\sum_{n\geq 0} a_n \Bigg(p_0 \left(4 \left(H_n^{(3)}\right){}^2-68 H_{2 n}^{(3)} H_n^{(3)}+289 \left(H_{2 n}^{(3)}\right){}^2+2 H_n^{(6)}-127 H_{2 n}^{(6)}\right) \\
		&\qquad +\frac{108 H_n^{(3)}}{(2 n+1)^2}-\frac{918 H_{2 n}^{(3)}}{(2 n+1)^2}-\frac{54}{(2 n+1)^5} \Bigg) \\
		&\qquad \stackrel{?}{=} 3072 \Im \Li_{4,1}(i,1)+1024 \Im \Li_{4,1}(i,-1)+2048 L_{-4}(4) \log (2) +\frac{3600 \zeta (3)^2}{\pi }-\frac{9434 \pi ^5}{945}.
	\end{aligned}
\end{gather}
\endgroup
In the last equation, we use the following notation for the multiple polylogarithm:
$$\Li_{s_1,s_2}(a,b) = \sum_{n_1>n_2\geq 1} \frac{a^{n_1}b^{n_2}}{n_1^{s_1}n_2^{s_2}}.$$ 
The first equation above is a special case of Sun's conjecture \cite[Conjecture~29]{sun2022conjectures} on Ramanujan's $1/\pi$ formula, which was confirmed by Zhou \cite[Section~3]{zhou2025notes} via modular equations; see also \cite{hou2024ramanujan}. Equations (4.22) and (4.23) were proved in \cite{li2024series}; (4.28) would follow from two conjectures of Sun in \cite{sun2023new, sun2022conjectures}; and (4.30) was also conjectured by Sun \cite[(5.14)]{sun2023new}. Such sums were also investigated in \cite{li2024series}, but only those that follow from $S_2$, thereby excluding most of the identities above. Despite using multiple WZ-seeds, the conjecture \cite[Conjecture~38]{sun2022conjectures}
$$\sum_{n\geq 0}  \left(\frac{1}{2^6}\right)^n \frac{(\frac12)_n^3}{(1)_n^3} (42n+5) \left(H_{2n}^{(3)} - \frac{43}{352} H_n^{(3)}\right) =\frac{555\zeta(3)}{77\pi} - \frac{32G}{11}$$ does not follow from the considerations above. This has been recently established in \cite{sun2026series} via Eichler integral of modular forms. \par 

A step-by-step computational guide for this example is available in the Mathematica package mentioned in the introduction. Readers who wish to see more formulas can consult the package. 
\end{romexample}

\begin{romexample}\label{Ex_8}
	To discover a harmonic extension of Ramanujan's $1/\pi$ series \cite{cohen2021rational, chudnovsky1988approximations},
	$$\sum_{n\geq 0} \left(\frac{1}{9}\right)^n \frac{(\frac12)_n (\frac14)_n (\frac34)_n}{(1)_n^3} (1+8n) = \frac{2\sqrt{3}}{\pi},$$
	we consider a WZ-seed that proves it, $$F(n,k) = \textsf{Seed10}\left(a+3 b-n+\frac{1}{2},2 b+1,c+k+n\right).$$
	From its WZ-mate, we obtain
	$$S:= S(a,b,c;n) := \frac{3^{-2 n} (2 a-6 b-6 c-8 n-1) \left(-a-3 b+\frac{1}{2}\right)_n \left(-a+3 b+\frac{1}{2}\right)_n \left(a+3 b+3 c+\frac{1}{2}\right)_{2 n}}{(1-2 a)_{2 n} (c+1)_n (2 b+c+1)_n}.$$
	\begin{conjecture}
		The generating function $$\sum_{N\geq 0} \dim V_N(S) t^N \stackrel{?}{=} \frac{1}{(1-t)^2(1-t^2)} = 1+2t+4t^2+6t^3+9t^4+12t^5+\cdots$$
	\end{conjecture}
	Set $$a_n = \left(\frac{1}{9}\right)^n \frac{(\frac12)_n (\frac14)_n (\frac34)_n}{(1)_n^3},\qquad p_0 = 1+8n.$$ We list some examples of identities in $V_N(S)$, including only the simplest-looking ones.
	\begin{align*}
		&\color{NavyBlue}\sum_{n\geq 0}a_n p_0 \left(1+\frac{1}{6} \left(-4 H_n+3 H_{2 n}\right)\right) \stackrel{?}{=} \frac{8 \log (2)}{\sqrt{3} \pi },\\
		&\color{NavyBlue}\sum_{n\geq 0}a_n \left(2+p_0(H_{4n}-H_n)\right)  = \frac{4 \sqrt{3} \log (2)}{\pi }+\frac{\sqrt{3} \log (3)}{\pi} \\
		&\color{NavyBlue}\sum_{n\geq 0}a_n p_0 \left(-5 H_n^{(2)}+18 H_{2 n}^{(2)}\right)  =\frac{\pi }{\sqrt{3}} \\
		&\sum_{n\geq 0}a_n \left(-4 H_n+3 H_{2 n}+\frac{1}{24} \left(32 (H_n)^2-48 H_n H_{2 n}+18 (H_{2 n})^2+5 H_n^{(2)}\right) p_0\right)  \stackrel{?}{=}  \frac{32 \log ^2(2)}{\sqrt{3} \pi }-\frac{37 \pi }{24 \sqrt{3}}, \\
		&\sum_{n\geq 0}a_n \left(-H_n+H_{4 n}+\frac{1}{16} \left(4 \left(H_n\right){}^2-8 H_n H_{4 n}+4 \left(H_{4 n}\right){}^2+H_n^{(2)}-4 H_{4 n}^{(2)}\right) p_0\right) \\ &\quad  \stackrel{?}{=}  -\frac{3 \sqrt{3} \pi }{16}+\frac{2 \sqrt{3} \log ^2(2)}{\pi }+\frac{\sqrt{3} \log ^2(3)}{8 \pi }+\frac{\sqrt{3} \log (3) \log (2)}{\pi } \\
		&\sum_{n\geq 0}a_n \left(-5 H_n^{(2)}+18 H_{2 n}^{(2)}+\frac{1}{6} \left(20 H_n H_n^{(2)}-15 H_{2 n} H_n^{(2)}+14 H_n^{(3)}-72 H_n H_{2 n}^{(2)}+54 H_{2 n} H_{2 n}^{(2)}-108 H_{2 n}^{(3)}\right) p_0 \right)  \\&\quad  \stackrel{?}{=}  \frac{45 L_{-3}(2)}{2}-\frac{94 \zeta (3)}{\sqrt{3} \pi }+\frac{4 \pi  \log (2)}{3 \sqrt{3}}, \\
		&\sum_{n\geq 0}a_n p_0 \left(50 (H_n^{(2)})^2-360 H_{2 n}^{(2)} H_n^{(2)}+648 (H_{2 n}^{(2)})^2+41 H_n^{(4)}-648 H_{2 n}^{(4)} \right)  \stackrel{?}{=}  -\frac{2 \pi ^3}{15 \sqrt{3}}.
	\end{align*}
\end{romexample}

\begin{remark}
(a) The first identity was conjectured by Sun \cite[Conjecture~14]{sun2022conjectures}; the third was proved in \cite{wei2024some}. \\
(b) The second identity also follows from a general conjecture of Sun \cite[Conjecture~29]{sun2022conjectures} on Ramanujan's $1/\pi$ formula, which was confirmed by Zhou \cite[Section~3]{zhou2025notes} via modular equations. The same applies to the second identity in the next example.
\end{remark}

\begin{romexample}\label{Ex_9}
		To discover a harmonic extension of Ramanujan's $1/\pi$ series \cite{cohen2021rational, chudnovsky1988approximations},
	$$\sum_{n\geq 0} \left(-\frac{3^3}{2^9}\right)^n \frac{(\frac12)_n (\frac16)_n (\frac56)_n}{(1)_n^3} (15+154n) = \frac{32\sqrt{2}}{\pi}$$ 
	we consider a WZ-seed that proves it, 	$$F(n,k) = \textsf{Seed1}\left(-a+b-c-n+\frac{1}{2},b-c+1,c+k+2 n\right).$$
	From its WZ-mate, we obtain
	$$S(a,b,c;n) := \frac{(-1)^n 2^{-5 n} P(n) (a+b-c+\frac{1}{2})_n (a-b+c+\frac{1}{2})_n (-a+b+c+\frac{1}{2})_{3 n}}{(a+1)_n (b+1)_{2 n+1} (c+1)_{2 n+1}},$$
	where
	\begin{multline*}
		P(n) = 15+ a^3+a^2 (-8 b-8 c+4)+8 b^3+12 b^2 (2 c+3)+b \left(24 c^2+72 c+46\right)+a \left(-8 b^2+24 b (2 c+1)-8 c^2+24 c+30\right) \\ +n \left(-8 a^2+8 a (6 b+6 c+13)+88 b^2+24 b (10 c+13)+88 c^2+312 c+214\right) \\ + n^2 (88 a+440 b+440 c+676)+8 c^3+36 c^2+46 c+616 n^3.
	\end{multline*}
	We conjecture that the generating function of $\dim V_N(S)$ is the same as in the previous example. We list some examples of identities in $V_N(S)$, including only the simplest-looking ones. To abbreviate the notation, set $$a_n =  \left(-\frac{3^3}{2^9}\right)^n \frac{(\frac12)_n (\frac16)_n (\frac56)_n}{(1)_n^3},\quad p_0 = 15+154 n,\quad p_1 = 61+110 n,\quad p_2 = 1+6 n,\quad T_n^{(s)} = 2^s H_{2 n}^{(s)}-H_n^{(s)}.$$
	\begin{align*}&\color{NavyBlue}\sum_{n\geq 0} a_n \left(\left(-H_n+H_{2 n}\right) p_0+\frac{p_1}{3 (1+2 n)} \right)  \stackrel{?}{=}  \frac{64 \sqrt{2} \log (2)}{\pi } \\
		&\color{NavyBlue}\sum_{n\geq 0} a_n \left(\frac{1}{3}+\frac{1}{154} \left(-H_n+T_{3 n}^{(1)}\right) p_0\right) = \frac{80 \sqrt{2} \log (2)}{77 \pi } \\
		&\color{NavyBlue}\sum_{n\geq 0} a_n \left(\frac{1}{15} \left(-4 H_n^{(2)}+15 H_{2 n}^{(2)}\right) p_0-\frac{p_2}{(1+2 n)^2}\right)  \stackrel{?}{=}  -\frac{16 \sqrt{2} \pi }{45} \\
		&\sum_{n\geq 0} a_n \left(\frac{1}{18} \left(9 H_n^2-18 H_n H_{2 n}+9 (H_{2 n})^2+H_n^{(2)}\right) p_0-\frac{\left(H_n-H_{2 n}\right) p_1}{3 (1+2 n)}-\frac{p_2}{2 (1+2 n)^2} \right)  \stackrel{?}{=}  \frac{64 \sqrt{2} \log ^2(2)}{\pi }-\frac{88 \sqrt{2} \pi }{27} \\
		&\sum_{n\geq 0} a_n \left(-2 \left(H_n-T_{3 n}^{(1)}\right)+\frac{1}{154} \left(3 (H_n)^2+H_n^{(2)}-6 H_n T_{3 n}^{(1)}+3 (T_{3 n}^{(1)})^2-3 T_{3 n}^{(2)}\right) p_0 \right)  \stackrel{?}{=}  \frac{1200 \sqrt{2} \log ^2(2)}{77 \pi }-\frac{16 \sqrt{2} \pi }{21} \\
		&\sum_{n\geq 0} a_n \Bigg(\frac{1}{3} \left(77 H_{2 n}^{(2)}+\frac{6 n-1}{(2 n+1)^3}\right) +\frac{\left(H_n-T_{3 n}^{(1)}\right) p_2}{2 (1+2 n)^2} -\frac{308 H_n^{(2)}}{45} \\&\quad  +\frac{1}{90} \left(12 H_n H_n^{(2)}+8 H_n^{(3)}-45 H_n H_{2 n}^{(2)}-60 H_{2 n}^{(3)}-12 H_n^{(2)} T_{3 n}^{(1)}+45 H_{2 n}^{(2)} T_{3 n}^{(1)}\right) p_0  \Bigg)  \\&\quad  \stackrel{?}{=}  \frac{512 L_{-8}(2)}{45}-\frac{832 \sqrt{2} \zeta (3)}{45 \pi }-\frac{8}{9} \sqrt{2} \pi  \log (2)
	\end{align*}
\end{romexample}

\begin{romexample}\label{Ex_10}
	To discover a harmonic extension of the remarkable $\zeta(5)$ formula (conjectured by Zhao\footnote{\url{https://mathoverflow.net/questions/281009}} and proved by the author \cite{au2025wilf}) \begin{equation}\label{aux_5}\sum_{n\geq 1} \left(-\frac{2^{10}}{5^5}\right)^n \frac{(1)_n^9}{(\frac12)_n^5 (\frac15)_n(\frac25)_n(\frac35)_n(\frac45)_n} \frac{30-425 n+2275 n^2-5600 n^3+5532 n^4 }{n^9} = -380928 \zeta(5),\end{equation} we consider the same WZ-seed that was used to prove it:
	$$F(n,k) = \textsf{Seed7}\left(a+n+\frac{1}{2},b-n,c+n+\frac{1}{2},d+n+\frac{1}{2},e+k+n+\frac{1}{2}\right).$$
	We focus on the term $G(n-1,0)$ in the formula $$\sum_{k\geq 0}F(0,k) + \sum_{n\geq 0} g(n) = \sum_{n\geq 1} G(n-1,0).$$
	Again, let $S(a,b,c,d,e;n)$ denote $G(n-1,0)$ with the factors independent of $n$ omitted. Its degree-zero coefficient is exactly the summand in \eqref{aux_5}. 	\begin{conjecture}
		The generating function
		$$\sum_{N\geq 0} \dim V_N(S) t^N \stackrel{?}{=} \frac{1}{(1-t)(1-t^2)(1-t^3)(1-t^4)(1-t^5)} = 1+t+2t^2+3t^3+5t^4+7t^5+10t^6+\cdots.$$
	\end{conjecture}
	We list some examples of identities in $V_N(S)$, including only the simplest-looking ones. To shorten the lengthy formulas, write
	$$\begin{aligned}a_n &= \left(-\frac{2^{10}}{5^5}\right)^n \frac{(1)_n^9}{(\frac12)_n^5 (\frac15)_n(\frac25)_n(\frac35)_n(\frac45)_n} \\
		p_0 &=30-425 n+2275 n^2-5600 n^3+5532 n^4 \\
		p_1 &= 54-680 n+3185 n^2-6720 n^3+5532 n^4\\
		p_2&=  510-7225 n+37835 n^2-87620 n^3+76632 n^4.
	\end{aligned}$$
	Then
	{\begin{align*}&\color{NavyBlue}\sum_{n\geq 1} a_n \left(\frac{p_1}{n^{10}} + \frac{p_0 \left(-3 H_n+2 H_{2 n}+H_{5 n}\right)}{n^9}\right) \stackrel{?}{=} 786432 \zeta(\overline{5},1)-294912 \zeta (3)^2-1523712 \zeta (5) \log (2)+\frac{32768 \pi ^6}{105} \\
		&\color{NavyBlue}\sum_{n\geq 1} a_n \left(\frac{p_0 \left(11 H_{2 n}^{(2)}-7 H_n^{(2)}\right)}{n^9}+\frac{p_2}{4 n^{11}}\right) \stackrel{?}{=} -3901440 \zeta (7) \\
		&\color{NavyBlue}\sum_{n\geq 1} a_n \left(\frac{p_0 \left(H_{2 n}^{(3)}-H_n^{(3)}\right)}{n^9}+\frac{238512 n^5-351796 n^4+213070 n^3-65625 n^2+10185 n-630}{24 n^{12} (2 n-1)}\right) \stackrel{?}{=} -\frac{512 \pi ^8}{9}\\
		&\color{NavyBlue}\sum_{n\geq 1} a_n \left(\frac{p_0 \left(13 H_n^{(4)}+49 H_{2 n}^{(4)}\right)}{n^9}+\frac{\splitfrac{-6044128 n^6+11644928 n^5-9540232 n^4}{+4214800 n^3-1052415 n^2+140065 n-7710}}{16 (1-2 n)^2 n^{13}}\right) \stackrel{?}{=}-14651392 \zeta (9)\\
		&\sum_{n\geq 1} a_n \Bigg(\frac{p_2 \left(11 H_{2 n}^{(2)}-7 H_n^{(2)}\right)}{10 n^{11}}+\frac{p_0 \left(2401 \left(H_n^{(2)}\right){}^2-7546 H_{2 n}^{(2)} H_n^{(2)}+5929 \left(H_{2 n}^{(2)}\right)^2-1542 H_n^{(4)}\right)}{245 n^9}+ \\ &\quad \frac{\splitfrac{153068992 n^6-315048448 n^5+270889824 n^4-123727032 n^3}{+31524339 n^2-4232797 n+232998}}{784 (1-2 n)^2 n^{13}}\Bigg) \stackrel{?}{=} -\frac{44552192}{5} \zeta (9)  \\
		&\color{NavyBlue}\sum_{n\geq 1} a_n \Bigg(\frac{p_0 \left(7 H_n^{(5)}+117 H_{2 n}^{(5)}\right)}{n^9}+\frac{\splitfrac{-43911744 n^7+112688016 n^6-121035008 n^5+71405368 n^4}{-25155570 n^3+5304255 n^2-618915 n+30690}}{96 n^{14} (2 n-1)^3}\Bigg)  \\ &\quad \color{NavyBlue}\stackrel{?}{=} -23617536 \zeta (5)^2-\frac{2048 \pi ^{10}}{9}
	\end{align*} }
\end{romexample}

	\bibliographystyle{plain} 
	\bibliography{../ref.bib} 
	
\end{document}